\documentclass[a4paper,12pt]{article}

\usepackage{amsmath, amssymb, amsthm}
\usepackage[dvips,final]{graphicx}
\usepackage[dvips,final]{epsfig}

\usepackage{xcolor}

\usepackage{titlesec}
\titleformat{\section}{\bfseries}{\thesection}{1em}{}
\titleformat{\subsection}{\itshape}{\thesubsection}{1em}{}

\newfont{\ctv}{msam10}

\newcommand{\bbox}{\mbox{\ctv \symbol{4}}}
\def\QED{{${}\hfill\bbox$}}
\newenvironment{pf}[1]{\par\vskip1mm{\noindent\it #1.}\ }{\QED\par
\vskip2mm}

\def\bpf{\begin{pf}}
\def\epf{\end{pf}}

\def\expe{\hbox{\rm e}}
\def\ve{\varepsilon}
\def\vp{\varphi}
\def\dd{\,\mathrm{d}}

\def\sign{\mathrm{\,sign}}
\def\supess{\mathop{\mathrm{\,sup\,ess}}}
\def\for{\mathrm{\ for\ }}
\def\ale{\mathrm{\ a.\,e.}}
\def\meas{\mathrm{meas\,}}
\def\play{\mathfrak{p}}
\def\on{^{(n)}}

\def\sumj{\sum_{i=1}^{n}}
\def\sumi{\sum_{i=0}^{n-1}}

\def\real{\mathbb{R}}
\def\nat{\mathbb{N}}

\def\io{\int_{\Omega}}
\def\ipo{\int_{\partial\Omega}}

\def\be{\begin{equation}\label}
\def\ee{\end{equation}}

\def\ber{\begin{eqnarray}}
\def\eer{\end{eqnarray}}

\def\bers{\begin{eqnarray*}}
\def\eers{\end{eqnarray*}}

\def\bpf{\begin{pf}}
\def\epf{\end{pf}}

\newtheorem{theorem}{Theorem}[section]
\newtheorem{lemma}[theorem]{Lemma}

\newtheorem{hypothesis}[theorem]{Hypothesis}
\newtheorem{proposition}[theorem]{Proposition}
\newtheorem{definition}[theorem]{Definition}

\begin{document}

\title{Degenerate diffusion with Preisach hysteresis\thanks{The first author is supported by the Austrian Science Fund (FWF) grants V662, Y1292, and F65, and by the OeAD WTZ grants CZ02/2022 and CZ09/2023; the second author is supported by the GA\v CR Grant No.~20-14736S, and by the European Regional Development Fund Project No.~CZ.02.1.01/0.0/0.0/16{\_}019/0000778.}}

\author{Chiara Gavioli
\thanks{Institute of Analysis and Scientific Computing, TU Wien, Wiedner Hauptstra\ss e 8-10, A-1040 Vienna, Austria, E-mail: {\tt chiara.gavioli@tuwien.ac.at}.}
\and Pavel Krej\v c\'{\i}
\thanks{Faculty of Civil Engineering, Czech Technical University, Th\'akurova 7, CZ-16629 Praha 6, Czech Republic, E-mail: {\tt Pavel.Krejci@cvut.cz}.}
}

\date{}

\maketitle

\begin{abstract}
Fluid diffusion in unsaturated porous media manifests strong hysteresis effects due to surface tension on the liquid-gas interface. We describe hysteresis in the pressure-saturation relation by means of the Preisach operator, which makes the resulting evolutionary PDE strongly degenerate. We prove the existence and uniqueness of a strong global solution in arbitrary space dimension using a special weak convexity concept.

\bigskip


\medskip


\end{abstract}


\section*{Introduction}\label{sec:intr}

This paper is a contribution to the discussion about the role of convexity in PDEs with hysteresis. We show here a typical situation, where even an approximate convexity concept can be sufficient for controlling the global evolution of a degenerate system in space and time.

\begin{figure}[htb]
\begin{center}
\includegraphics[width=9cm]{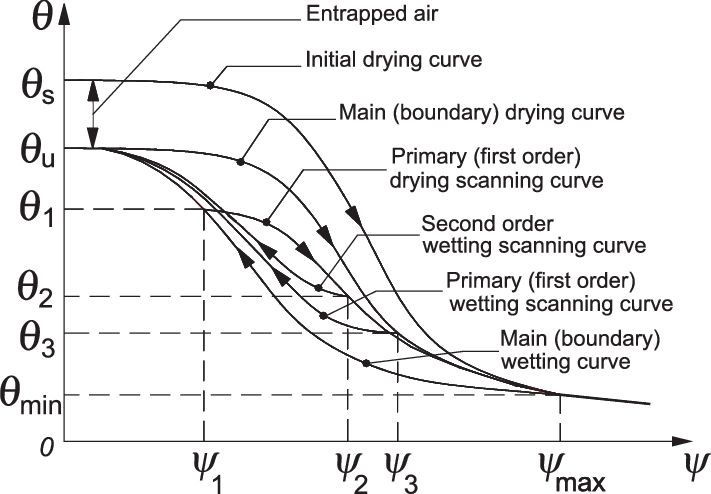}

\caption{Typical experimental hysteresis dependence in porous media between the logarithm soil suction $\psi$ and the volumetric water content $\theta$.}
\label{f1}
\end{center}
\end{figure}

As a model situation, we consider the problem of fluid diffusion in porous media governed by the PDE
\be{e0}
s_t - \Delta u = h(x,t)
\ee
with space variable $x \in \Omega$, where $\Omega \subset \real^N$ for $N \in \nat$ is a bounded open domain, and time variable $t \in (0,T)$ for some $T > 0$. The unknown function $u = u(x,t) \in \real$ represents the pressure, $s = s(x,t) \in (0,1)$ is the relative saturation of the fluid in the pores, $\Delta$ is the Laplacian in $x$, and $h(x,t) \in \real$ is the source density. As a result of surface tension on the air-liquid interface and adhesion on the liquid-solid contact surface, the pressure-saturation relation exhibits rate-independent hysteresis, see Figure~\ref{f1} taken from \cite{pfb}, where $\theta$ represents the volumetric water content (that is, saturation), and $\psi$ is the logarithm soil suction, which can be interpreted as a decreasing function of the pressure. We follow the modeling hypotheses of \cite{alb,bsch} and represent hysteresis by a \textit{Preisach operator} $G$ in the form
\be{e0a}
s(x,t) = G[u](x,t).
\ee
System \eqref{e0}--\eqref{e0a} is considered with
a given right-hand side $h(x,t)$, with a given initial condition
\be{e2}
u(x,0) = u^0(x),
\ee
and with boundary condition
\be{e1}
-\nabla u\cdot n = b(x) (u - u^*(x,t))
\ee
on $\partial\Omega$ with a given boundary source $u^*(x,t)$ and given boundary permeability $b(x)\ge 0$.

The Preisach operator will be rigorously defined below in Definition~\ref{dpr}. Let us only mention at this point that
the time evolution described by equation \eqref{e0} is doubly degenerate: in typical situations, the function $s(x,t) = G[u](x,t)$  is bounded independently of the evolution of $u$, so that no a~priori lower bound for $u_t$ is immediately available. Furthermore, at every point $x \in \Omega$ and every time $t_0$ where $u_t$ changes sign (which the engineers call a {\em turning point\/}), we have
\be{e6a}
u_t(x,t_0-\delta)\cdot u_t(x,t_0+\delta)<0 \ \ \forall \delta\in (0,\delta_0(x)) \ \Longrightarrow \ \liminf_{\delta \to 0+} \frac{G[u]_t(x,t_0+\delta)}{u_t(x,t_0+\delta)} = 0,
\ee
see Figure~\ref{f1}, so that the knowledge of $G[u]_t$ alone does not give a complete information about $u_t$. In order to control the time derivatives of $u$, we have to restrict our considerations to the so-called convexifiable Preisach operators, and refer to deeper properties of the \textit{memory} of the system: the concept of hysteresis memory is related to the fact that the instantaneous output value may depend not only on the instantaneous input value and the initial condition, but also on other input values in the history of the process. For the Preisach model, memory is completely characterized by a real function of one variable called the memory curve specified below after Definition~\ref{dpr}. In the case under consideration, if for example $\Delta u^0(x) + h(x,0) > 0$ and the initial memory (that is, the memory created prior to $t=0$) corresponds to a turning point, then the process described by \eqref{e0} cannot start at all at $t=0$, and no solution can be expected to exist. For this reason, the initial memory curve at $t=0$ has to satisfy a suitable memory compatibility condition stated below in Hypothesis~\ref{hym} and exploited in Subsection~\ref{sec:init}.

The highly degenerate problem we consider here does not fall within the classical theory of parabolic equations with hysteresis developed by Visintin, see \cite{vis}, where the left-hand side of \eqref{e0} contains an additional term $c u_t$ with $c > 0$. In such a case, an immediate estimate for $u_t$ is available from the very beginning, and a lot of existence results even for more general (nonlinear) elliptic operators and boundary conditions have been proved.

The structure of the paper is the following. In Section~\ref{sec:stat} we list the definitions of the main concepts including convexifiable Preisach operators, and state the main existence and uniqueness Theorem~\ref{t1}. In Section~\ref{sec:disc} we state and solve the time-discrete counterpart of \eqref{e0a}--\eqref{e1}. Section~\ref{sec:itau} is devoted to estimates for solutions to the time-discrete problem which are independent of the time step. Finally, in Section~\ref{sec:proof} we pass to the limit as the time step tends to zero and check that the limit is the unique solution to the original problem.


\section{Statement of the problem}\label{sec:stat}

\begin{definition}\label{dpr}
Let $\lambda \in L^\infty(\Omega \times (0,\infty))$ be a given function which we call {\em initial memory curve\/} and which has the following properties: 
\begin{align}\label{e6}
&|\lambda(x, r_1) - \lambda(x, r_2)| \le |r_1 - r_2| \ \ale \ \forall r_1, r_2 \in (0,\infty),\\ \label{e6b}
&\exists \Lambda > 0:\ \lambda(x,r) = 0\ \for r\ge \Lambda \ \mbox{ and a.\,e. } x \in \Omega.
\end{align}
For a given $r>0$, we call the {\em play operator with threshold $r$ and initial memory $\lambda$} the mapping which with a given function $u \in L^2(\Omega; W^{1,1}(0,T))$ associates the solution $\xi^r\in L^2(\Omega; W^{1,1}(0,T))$ of the variational inequality
\be{e4a}
|u(x,t) - \xi^r(x,t)| \le r, \quad \xi^r_t(x,t)(u(x,t) - \xi^r(x,t) - z) \ge 0 \ \ale \ \ \forall z \in [-r,r],
\ee
with a given initial memory curve
\be{e5}
\xi^r(x,0) = \lambda(x,r)
\ee
with $\lambda(x,0) = u(x,0)$ for a.\,e. $x \in \Omega$, and we denote
\be{e4}
\xi^r(x,t) = \play_r[\lambda,u](x,t).
\ee
Given a measurable function $\rho :\Omega\times(0,\infty)\times \real \to [0,\infty)$ and a constant $\bar G \in \real$, the Preisach operator $G$ is defined as a mapping $G: L^2(\Omega; W^{1,1}(0,T))\to L^2(\Omega; W^{1,1}(0,T))$ by the formula
\be{e3}
G[u](x,t) = \bar G + \ \int_0^\infty\int_0^{\xi^r(x,t)} \rho(x,r,v)\dd v\dd r.
\ee
The Preisach operator is said to be {\em regular\/} if the density function $\rho$ of $G$ in \eqref{e3} belongs to $L^1(\Omega\times (0,\infty )\times \real)\cap L^\infty(\Omega\times (0,\infty )\times \real)$, and there exist a constant $\rho_1$ and a decreasing function $\rho_0: [0,\infty) \to (0,\infty)$ such that for every $U>0$ we have
\be{e3a}
\rho_1 > \rho(x,r,v)> \rho_0(U) > 0 \ \for \ale\ (x,r,v)\in \Omega\times (0,U) \times (-U,U).
\ee
\end{definition}

The original construction of the Preisach operator, which goes back to \cite{prei}, see also \cite{bert,bsch,vis}, represents hysteresis by a continuous distribution of relays, which can switch between $-1$ and $+1$. The equivalence of the two definitions is proved in \cite[Lemma~II.3.6]{book}. The advantage of the play operator representation is that analytical properties (Lipschitz continuity with respect to $\sup$-norm and $BV$-norm, energy inequalities, etc.) of the Preisach operator can easily be deduced from the corresponding properties of the play operator.

\begin{figure}[htb]
\begin{center}
\includegraphics[width=16cm]{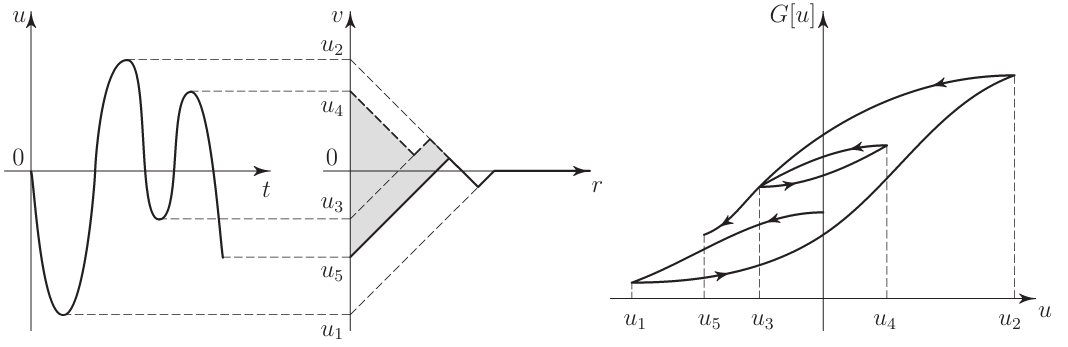}
\caption{Evolution of the memory curve and of the hysteresis diagram associated with the input sequence $0 {\to} u_1 {\to} u_2 {\to} u_3 {\to} u_4 {\to} u_5$.}
\label{f2}\end{center}
\end{figure}

In Definition~\ref{dpr} the variable $r$ is the \textit{memory variable}, and at each point $x\in\Omega$ and time $t\ge 0$, the so-called \textit{memory curve} $r \mapsto \play_r[\lambda,u](x,t)$ describes the whole memory created during the interval $[0,t]$. The left picture in Figure~\ref{f2} illustrates how the memory curve in the $(r,v)$-plane is updated when $\lambda(x,r) \equiv 0$, and the input $u$ decreases from $u_0 = 0$ to some value $u_1 < 0$, then it increases from $u_1$ to $u_2$, then decreases from $u_2$ to $u_3$, then increases from $u_3$ to $u_4$, and then decreases again from $u_4$ to $u_5$. The dashed lines in the $(r,v)$-plane correspond to parts of the memory which have been erased during the evolution. This is typical of \textit{minor loops}, where the process returns back to the starting point and continues as if no turn had taken place, see the transition $u_4 \to u_5$. The memory curve is typically piecewise linear with slopes $+1, -1, 0$. If we call `corners' the points on the memory curve where the slope changes, we see that during a continuous monotone evolution of $u$ only the left corner moves to the right in a possibly discontinuous way. The $r$-coordinate of this active corner is called the {\em memory depth\/} or {\em active memory level\/}. Geometrically, the instantaneous value of $G[u](x,t) - \bar G$ at a point $x \in \Omega$ and time $t$ is represented by the integral of  $\rho(x,r,v)$ over the 2D domain of parameters $(r,v)$ between the $r$-axis and the memory curve $v = \play_r[\lambda,u](x,t)$, cf.~\eqref{e3}. For more details, see \cite[Sections~II.2--II.3]{book}.

The rate-independence of $G$ implies that for input functions $u(x,t)$ which are monotone in a time interval $t\in [a,b]$, a regular Preisach operator $G$ can be represented by a superposition operator
\begin{equation}\label{GB}
G[u](x,t) = G[u](x,a) + B(x, u(x,t)) \ \mbox{ for } t \in [a,b]
\end{equation}
with an increasing function $u \mapsto B(x,u)$ called a {\em Preisach branch\/}. Indeed, the branches may be different at different points $x$ and different intervals $[a,b]$. The branches corresponding to increasing inputs are said to be {\em ascending\/}, the branches corresponding to decreasing inputs are said to be {\em descending\/}. Referring again to the geometric interpretation mentioned above, the branch is represented by the integral of $\rho(x,r,v)$ over the area between the memory curves $v = \play_r[\lambda,u](x,a)$ and $v = \play_r[\lambda,u](x,t)$. For example, the area of the moving gray zone in Figure~\ref{f2} measured in terms of the density $\rho$ describes the evolution of the descending branch from $u_4$ to $u_5$. For a regular Preisach operator, the functions $u \mapsto B(x,u)$ are always continuous and their derivatives $u \mapsto B_u(x,u)$ with respect to $u$ have locally bounded variation, see \cite[Lemma~II.3.18 and Proposition~II.4.14]{book}. In the situation of Figure~\ref{f2}, when $u$ crosses the value $u_3$ on the way from $u_4$ to $u_5$, the memory created by the local maximum at $u_3$ is erased, and $B_u(x,u)$ has a downward jump at $u = u_3$. In general, the derivative of an ascending branch can only jump up, and the derivative of a descending branch can only jump down. In addition, for regular Preisach operators, the locally positive lower bound for the density $\rho$ guarantees at least local convexity/concavity of the branches, see again \cite[Sections~II.3 and II.4]{book}. It is also easy to understand why at a turning point $u_*$ the right derivative $B_u(x,u_*+)$ vanishes, and \eqref{e6a} holds. Indeed, if $|u - u_*| = \ve$, then the area of the triangle between the memory curves is $\ve^2/4$. 

As anticipated in the Introduction, in order to control the time derivatives of $u$ we have to restrict our analysis to convexifiable Preisach operators. The precise definition is given below.

\begin{definition}\label{dpc}
A regular Preisach operator is said to be {\em uniformly counterclockwise convex\/} on an interval $[-U,U]$ for $U>0$ if for all inputs $u:\Omega\times [0,T] \to [-U,U]$ all ascending branches are uniformly convex and all descending branches are uniformly concave.
\\[1mm]
A regular Preisach operator $G$ is called {\em convexifiable on\/} $[-U,U]$ if there exist a uniformly counterclockwise convex Preisach operator $P$ on $[-U,U]$, positive constants $g_*(U), g^*(U), \bar g(U)$, and a twice continuously differentiable mapping $g$ of $[-U,U]$ onto $[-U, U]$ such that
\be{hg}
g(0) = 0,\quad 0 <g_*(U) \le g'(u) \le g^*(U), \quad |g''(u)|\le \bar g(U)\ \  \forall u\in (-U,U),
\ee
and $G[u]= (P\circ g)[u] := P[g(u)]$ for all inputs $u:\Omega\times [0,T] \to [-U,U]$.\\[1mm]
A regular Preisach operator is called {\em convexifiable\/} if it is convexifiable on $[-U,U]$ for every $U>0$.
\end{definition}

A typical example of a uniformly counterclockwise convex operator on $[-U,U]$ for every $U>0$ is the so-called {\em Prandtl-Ishlinskii operator\/} characterized by positive density functions $\rho(x,r)$ independent of $v$, see \cite[Section~II.3]{book}. Operators of the form $P\circ g$ with a Prandtl-Ishlinskii operator $P$ and an increasing function $g$ are often used in control engineering because of their explicit inversion formulas, see \cite{al,viso,kk}. They are called the {\em generalized Prandtl-Ishlinskii operators\/} (GPI) and represent an important subclass of convexifiable Preisach operators.

Unlike in the GPI case, where the function $g$ is fixed, the functions $g$ in the definition of a convexifiable Preisach operator may be different in different intervals $[-U,U]$. Note also that for every regular Preisach operator $P$ and every Lipschitz continuous increasing function $g$, the superposition operator $G = P\circ g$ is also a regular Preisach operator, and there exists an explicit formula for its density, see \cite[Proposition~2.3]{error}.

We show here another class of convexifiable Preisach operators which are not GPI. A general characterization of convexifiable Preisach operators is, however, still missing.

\begin{proposition}\label{pcon}
Let $G$ be a regular Preisach operator according to Definition~\ref{dpr}, and let its density $\rho$ satisfy for a.\,e. $(x,r,v)\in \Omega\times (0,\infty)\times \real$ the condition
\be{rgr}
\rho_v(x,r,v) = -\phi(v) \rho(x,r,v)
\ee
with an odd nondecreasing continuous function $\phi: \real \to \real$. Then for every $U>0$ there exists an odd increasing twice continuously differentiable function $\hat g: [-U, U] \to [-U, U]$, $\hat g(0) = 0$, $\hat g(U) = U$, $\hat g$ convex on $[0,U]$, concave on $[-U,0]$, $\hat g'(u) \ge \hat g_*(U)$ for some $\hat g_*(U)>0$ and for all $u \in [-U,U]$, and such that the Preisach operator $P = G\circ \hat g$ is uniformly counterclockwise convex on $[-U, U]$.
\end{proposition}

Indeed, the convexifiability condition is then fulfilled for $g = \hat g^{-1}$. We shall see that identity \eqref{rgr} is crucial. Its meaning is that the decay of $\rho(r,v)$ to $0$ as $v \to \pm\infty$ is somehow controlled by a function $\phi$ independently of $r$. The condition is probably not optimal, but still a large class of Preisach operators of this form is used in engineering practice, where the $x$-dependence of $\rho$ is usually neglected. A typical example is $\rho(r,v) = \alpha (r) \expe^{-\beta(v)}$ with a continuous function $\alpha>0$ and an even convex function $\beta \ge 0$, which satisfies the assumptions with $\phi(v) = \beta'(v)$ (see, e.\,g., \cite{bert}). It is also easy to see that if $\tilde G$ is of the form \eqref{e3} with density $\tilde\rho$ sufficiently ``close'' to $\rho$ in suitable norm and if $G$ satisfies the hypotheses of Proposition~\ref{pcon}, then $\tilde G$ is convexifiable, too.

\bpf{Proof}
We omit the $x$-dependence which plays no role in the argument. We also give all the details of the proof only for the case of the major ascending branch of $P$. The argument for the other cases (descending branches, minor loops) is fully analogous and we leave the details to the reader. As it has been mentioned in the discussion before Definition~\ref{dpc} referring to \cite{book}, branches of minor loops of regular Preisach operators may exhibit jumps in derivatives which always point upward for ascending branches and downward for descending branches, so that convexity/concavity is preserved, as illustrated in Figure~\ref{f2}.

Let $U > 0$. We consider the major ascending branch of $P$ in $[-U, U]$ given by the formula, see \cite[Lemma II.3.19]{book},
\be{n0}
B(u) = \int_0^{({\hat g}(u) + U)/2}\int_{-U+r}^{{\hat g}(u) - r} \rho (r, v)\dd v\dd r, \quad u \in [-U, U],
\ee
with $\hat g : [-U,U] \to [-U,U]$ to be chosen in such a way that the derivative of $B$
\be{n1}
B'(u) = {\hat g}'(u) \int_0^{({\hat g}(u) + U)/2}\rho (r, {\hat g}(u) - r)\dd r
\ee
is increasing in $[-U,U]$. For the second derivative of $B$ we get
\begin{align}\nonumber 
B''(u) &= {\hat g}''(u) \int_0^{({\hat g}(u) + U)/2}\rho (r, {\hat g}(u) - r)\dd r\\ \label{n2}
&\quad  + ({\hat g}'(u))^2 \left(\int_0^{({\hat g}(u) + U)/2}\rho_v (r, {\hat g}(u) - r)\dd r + \rho\Big(\frac{{\hat g}(u) + U}{2}, \frac{{\hat g}(u) - U}{2}\Big) \right).
\end{align}
Since ${\hat g}$ has to be odd, we have to distinguish the cases $u>0$ and $u<0$. By \eqref{e3a} we have $\rho \ge \rho_0(U) >0$ as long as $|u|\le U$. Hence, the function $B$ is uniformly convex on $[0,U]$ provided $\hat g'(u) \ge \hat g_*(U)$ for some $\hat g_*(U) > 0$ to be identified and for all $u \in (0,U)$, and
\be{n3} 
{\hat g}''(u) \int_0^{({\hat g}(u) + U)/2}\rho (r, {\hat g}(u) - r)\dd r
+ ({\hat g}'(u))^2 \int_0^{({\hat g}(u) + U)/2}\rho_v (r, {\hat g}(u) - r)\dd r\ge 0
\ee
for $u \in (0,U)$. For $u<0$ we set ${\bar u} = -u$ and use the identities ${\hat g}(u) = -{\hat g}{(\bar u)}$, ${\hat g}'(u) = {\hat g}'{(\bar u)}$, ${\hat g}''(u) = -{\hat g}''{(\bar u)}$ to get from \eqref{n2} a condition similar to \eqref{n3}, namely
\be{n4} 
-{\hat g}''{(\bar u)} \int_0^{(U - {\hat g}{(\bar u)})/2}\rho (r, -{\hat g}{(\bar u)} - r)\dd r + ({\hat g}'{(\bar u)})^2 \int_0^{(U - {\hat g}{(\bar u)})/2}\rho_v (r, -{\hat g}{(\bar u)} - r)\dd r \ge 0
\ee
for ${\bar u} \in (0,U)$.

We define ${\hat g}: [0,U] \to [0,U]$ to be the solution of the differential equation
\be{eqdif}
{\hat g}''(u) = \phi({\hat g}(u)) ({\hat g}'(u))^2, \quad {\hat g}(0) = 0, \ {\hat g}(U) = U
\ee
and claim that the convexity of $B$ is guaranteed in the sense that both inequalities \eqref{n3}--\eqref{n4} are satisfied. It suffices to check that for all $u \in (0,U)$ we have simultaneously
\begin{align}\label{n5}
\int_0^{({\hat g}(u) + U)/2} \left(\phi({\hat g}(u))\rho (r, {\hat g}(u) - r) + \rho_v (r, {\hat g}(u) - r)\right)\dd r & \ge 0, \\ \label{n6}
\int_0^{(U- {\hat g}(u))/2} \left(- \phi({\hat g}(u))\rho (r, -{\hat g}(u) - r) + \rho_v (r, -{\hat g}(u) - r)\right)\dd r & \ge 0.
\end{align}
Note first that in \eqref{n5} we have $\rho_v (r, {\hat g}(u) - r) \ge 0$ for $r \ge {\hat g}(u)$, and from \eqref{rgr} it follows that
\begin{align*}
&\int_0^{{\hat g}(u)} \Big(\phi({\hat g}(u))\rho (r, {\hat g}(u) - r) + \rho_v (r, {\hat g}(u) - r)\Big)\dd r\\
&\quad  = \int_0^{{\hat g}(u)} \Big(\phi({\hat g}(u)) - \phi({\hat g}(u) - r)\Big)\rho (r, {\hat g}(u) - r) \dd r \ge 0,
\end{align*}
which implies \eqref{n5}. In \eqref{n6} we have by \eqref{rgr} that
$$
\rho_v (r, -{\hat g}(u) - r) = \phi({\hat g}(u) + r)\, \rho (r, -{\hat g}(u) - r),
$$
hence,
\begin{align*}
&\int_0^{(U- {\hat g}(u))/2} \Big(-\phi({\hat g}(u))\rho (r, -{\hat g}(u) - r) + \rho_v (r, -{\hat g}(u) - r)\Big)\dd r\\
&\quad  = \int_0^{(U- {\hat g}(u))/2} \Big(-\phi({\hat g}(u)) + \phi({\hat g}(u) + r)\Big)\rho (r, -{\hat g}(u) - r) \dd r \ge 0,
\end{align*}
so that \eqref{n6} holds, which completes the convexity statement.

It remains to evaluate the solution of the differential equation \eqref{eqdif}, which can be rewritten in the form
\be{eqd1}
\log {\hat g}'(u) - \Phi({\hat g}(u)) = \log C,
\ee
where $\Phi(v) = \int_0^v \phi(s)\dd s$, and $C>0$ is a constant which is to be identified from the boundary condition ${\hat g}(U)=U$. We rewrite \eqref{eqd1} in the form
\be{eqd2}
{\hat g}'(u) \expe^{-\Phi({\hat g}(u))} = C.
\ee
For $v\ge 0$ put $\hat\Phi(v) = \int_0^v \expe^{-\Phi(s)} \dd s$. Then from the boundary condition ${\hat g}(0) = 0$ we get $\hat\Phi({\hat g}(u)) = Cu$. The function $\hat\Phi$ is bounded because of the superlinear growth of $\Phi$, it is increasing, and $\hat\Phi(0) = 0$. Hence, the solution to \eqref{eqdif} can be expressed in the form
$$
{\hat g}(u) = \hat\Phi^{-1}(Cu) \ \for \ u\in [0,U], \quad C = \frac1U \hat\Phi(U).
$$
From \eqref{eqd2} if follows that $\hat g'(u) \ge C > 0$, hence choosing $0 < \hat g_*(U) \le C$ the proof of Proposition~\ref{pcon} is completed.
\epf

For the data we prescribe the following regularity.

\begin{hypothesis}\label{hyd}
We assume that $\Omega \subset \real^N$ is of class $C^{1,1}$, and that the functions $h,u^0,b,u^*$ in \eqref{e0} and \eqref{e2}--\eqref{e1} are such that
\begin{align*}
h &\in L^\infty(\Omega\times (0,T)) \cap W^{1,2}(0,T; L^2(\Omega)),\\
u^0 &\in W^{2,2}(\Omega),\ \Delta u^0 \in L^\infty(\Omega),\\
b & \in W^{1,\infty}(\partial\Omega), \ b(x) \ge 0 \, \ale, \ \ipo b(x) \dd s(x) > 0,\\
u^* &\in L^\infty(\partial\Omega\times (0,T)) \cap W^{1,2}(\partial\Omega\times (0,T)).
\end{align*}
\end{hypothesis}

As it has been mentioned in the Introduction, even a local solution to Problem \eqref{e0}--\eqref{e2} may fail to exist if for example $\lambda(x,r) 
\equiv 0$ and $\Delta u^0(x) + h(x,0) \ne 0$. Then $t=0$ is a turning point for all $x \in \Omega$ and there is no way to satisfy \eqref{e0} in any sense. We therefore need an initial memory compatibility condition which we state in the following way.

\begin{hypothesis}\label{hym}
There exists a constant $L>0$ and a function $r_0 \in L^\infty(\Omega)$ such that the following initial compatibility conditions hold:
\begin{align} \label{c0}
\lambda(x,0) &= u^0(x) \ \ale,\\ \label{c1}
\frac1L\sqrt{|\Delta u^0(x) + h(x,0)|} &\le r_0(x) \le \Lambda, \\ \label{c2}
-\frac{\partial}{\partial r} \lambda(x,r) &\in \sign(\Delta u^0(x) + h(x,0)) \ \ale\ \for r\in (0,r_0(x)), \\ \label{c2a}
-\nabla u^0(x)\cdot n &= b(x) (u^0(x) - u^*(x,0)) \ \ale\ \mbox{\em on }\, \partial\Omega.
\end{align}
\end{hypothesis}

For each $h$ and $u^*$, the class of admissible initial conditions $u^0$ in terms of Hypothesis~\ref{hym} contains for example all solutions of the equation $\Delta u^0 = h^*$ with $h^* \in L^\infty(\Omega)$ and with boundary condition \eqref{c2a}. The value $r_0$ in \eqref{c1}--\eqref{c2} represents the memory depth at time $t=0$ in agreement with the discussion after Definition~\ref{dpr}. The meaning of conditions \eqref{c0}--\eqref{c2} is that for large values of $\Delta u^0(x) + h(x,0)$, the initial memory $\lambda$ has to go deeper in the memory direction. We shall see in Subsection~\ref{sec:init} that Hypothesis~\ref{hym} guarantees the existence of some previous admissible history of the process prior to the time $t=0$ and ensures the existence of a continuation for $t>0$.

Problem \eqref{e0}--\eqref{e1} is to be understood in variational form
\be{e0v}
\io G[u]_t\vp \dd x + \io\nabla u\cdot\nabla \vp\dd x + \ipo b(x)(u-u^*)\vp \dd s(x) = \io h\,\vp \dd x
\ee
for every test function $\vp \in W^{1,2}(\Omega)$. The main objective of this paper is the proof of the following result.

\begin{theorem}\label{t1}
Let Hypotheses~\ref{hyd} and \ref{hym} be satisfied, and let $T>0$ be given. Let $G$ be a convexifiable Preisach operator according to Definition~\ref{dpc}.
Then Problem \eqref{e0v} admits a unique strong solution $u \in L^2(0,T; W^{2,2}(\Omega))$ such that $\nabla u \in L^\infty(0,T; L^2(\Omega; \real^N))$, $u_t \in L^q(\Omega\times (0,T))$ for all $1 \le q < 1+ (2/N)$, $G[u]_t \in L^2(\Omega\times (0,T))$.
\end{theorem}

Convexity as a geometric property of hysteresis operators allows to derive higher order a~priori estimates which can be sufficient for both hyperbolic and parabolic PDEs with hysteresis to prove the existence of solutions, see \cite{det,jana,mz}. The main novelty of Theorem~\ref{t1} is that here the process is allowed to take place in the full data range and no restriction to the convexity domain of $G$ is required. This is possible thanks to a uniform maximum principle in Proposition~\ref{p2} below which guarantees that solutions to the time-discrete approximations admit a uniform $L^\infty$-bound $\bar U$. By convexifiability hypothesis, we convexify the operator $G$ on $[-\bar U, \bar U]$ by means of a suitable function $g$, and apply the convexity argument to the operator $P = G\circ g^{-1}$.


\section{Discretization}\label{sec:disc}

For an input sequence $\{u_i: i \in \nat\cup\{0\}\} \subset L^\infty(\Omega)$ we define the time-discrete Preisach operator by a formula of the form \eqref{e3}, namely,
\be{de3}
G[u]_i(x) = \bar G + \  \int_0^{\infty}\int_0^{\xi^r_i(x)} \rho(x,r,v)\dd v\dd r
\ee
where $\xi^r_i$ denotes the output of the time-discrete play operator
\be{de4}
\xi^r_i(x) = \play_r[\lambda,u]_i(x)
\ee
defined as the solution operator of the variational inequality
\be{de4a}
|u_i(x) - \xi^r_i(x)| \le r, \quad (\xi^r_i(x) - \xi^r_{i-1}(x))(u_i(x) - \xi^r_i(x) - z) \ge 0 \quad \forall i\in \nat \ \ \forall z \in [-r,r],
\ee
with a given initial memory curve
\be{de5}
\xi^r_0(x) = \lambda(x,r), \quad \lambda(x,0) = u_0(x),
\ee
similarly as in \eqref{e4a}--\eqref{e5}. Note that the discrete variational inequality \eqref{de4a} can be interpreted as weak formulation of \eqref{e4a} for piecewise constant inputs in terms of the Kurzweil integral, and details can be found in \cite[Section 2]{ele}. This is, however, not the objective of this paper.

Choosing $z= u_{i-1}(x) - \xi^r_{i-1}(x)$ in \eqref{de4a}, we see that the mapping $u_i(x) \mapsto \xi^r_i(x)$ defined by \eqref{de4a} is monotone in the sense
$$
(\xi^r_i(x) - \xi^r_{i-1}(x))(u_i(x) - u_{i-1}(x)) \ge (\xi^r_i(x) - \xi^r_{i-1}(x))^2,
$$
hence $|\xi^r_i(x) - \xi^r_{i-1}(x)| \le |u_i(x) - u_{i-1}(x)|$. This and \eqref{e3a} yield
\be{de5a}
\frac1C|G[u]_i(x) - G[u]_{i-1}(x)|^2 \le (u_i(x) - u_{i-1}(x))(G[u]_i(x) - G[u]_{i-1}(x)) \le C|u_i(x) - u_{i-1}(x)|^2,
\ee
with a constant $C>0$ depending only on $\rho_1$.

We further have the energy inequality which can be stated as follows. For $\xi \in \real$ put
\be{psi}
\psi(x,r,\xi) = \int_0^\xi\rho(x,r,v)\dd v, \quad \Psi (x,r,\xi) = \int_0^\xi v \rho(x,r,v)\dd v.
\ee
Since $\psi$ is an increasing function of $\xi$, it follows from \eqref{de4a} for $z=0$ that
\begin{equation}\label{ene0a}
\big(\psi(x,r,\xi^r_i(x)) - \psi(x,r,\xi^r_{i-1}(x))\big) (u_i(x) - \xi^r_i(x)) \ge 0
\end{equation}
for all $i \in \nat$. In both cases $\xi^r_{i-1}(x)<\xi^r_{i}(x)$ and $\xi^r_{i-1}(x)\ge \xi^r_{i}(x)$ we have that
\be{psi2}
\Psi(x,r,\xi^r_{i}(x)) {-} \Psi(x,r,\xi^r_{i-1}(x)) = \int_{\xi^r_{i-1}(x)}^{\xi^r_{i}(x)} v \rho(x,r,v)\dd v \le
\xi^r_i(x)\big(\psi(x,r,\xi^r_i(x)) {-} \psi(x,r,\xi^r_{i-1}(x))\big),
\ee
hence, by \eqref{ene0a},
\begin{equation}\label{ene0b}
\big(\psi(x,r,\xi^r_i(x)) - \psi(x,r,\xi^r_{i-1}(x))\big) u_i(x) \ge \Psi(x,r,\xi^r_{i}(x)) - \Psi(x,r,\xi^r_{i-1}(x)).
\end{equation}
We now define the energy functional $E$ associated with $G$ by the formula
\be{ene2}
E[u]_i(x) = \int_0^{\infty}\Psi(x,r,\xi^r_i(x))\dd r = \int_0^{\infty}\int_0^{\xi^r_i(x)} v\rho(x,v,r)\dd v\dd r.
\ee
Integrating \eqref{ene0b} in $r$ from $0$ to $\infty$ and invoking \eqref{de3}, we obtain the energy inequality
\be{ene1}
(G[u]_i(x) - G[u]_{i-1}(x))\,u_i(x) \ge E[u]_i(x) - E[u]_{i-1}(x)
\ee
for all $i\in \nat$. Note that the definition of $\xi^r_0$ in \eqref{de5} and the assumptions on $\rho$ and $\lambda$ in Definition~\ref{dpr} imply that the initial energy $E[u]_0$ satisfies the uniform estimate
\be{iniene}
0 \le E[u]_0(x) = \int_0^\Lambda\int_0^{\lambda(x,r)}v \rho(x,v,r)\dd v\dd r \le \frac{\rho_1}{2}\int_0^\Lambda \lambda^2(x,r)\dd r \le C,
\ee
with a constant $C$ depending only on $\rho_1$ and $\Lambda$. We now replace \eqref{e0} with its time-discrete counterpart
\be{e0e}
\frac1\tau \big(G[u]_i - G[u]_{i-1}\big) - \Delta u_i = h_i(x)
\ee
with a fixed time step $\tau>0$ which will be specified below in \eqref{tz}. The unknown functions are $\{u_i: i = 1, \dots, {n}\}$ with a given initial condition $u_0 = u^0$ and boundary conditions
\be{de1}
-\nabla u_i\cdot n = b(x) (u_i - u^*_i(x)),
\ee
where the discrete data
\be{drh}
h_i(x) = h(x,t_i), \quad u^*_i(x) = u^*(x,t_i)
\ee
are associated with a division $\{t_i = i\tau: i = 0,1,\dots, {n}\}$, ${n}\tau = T$ of the interval $[0,T]$.  We state the problem in variational form
\be{de0v}
\io \left(\frac1\tau(G[u]_i - G[u]_{i-1})\vp + \nabla u_i\cdot\nabla\vp\right)\dd x + \ipo b(x)(u_i - u^*_i)\vp \dd s(x) = \io h_i\vp\dd x
\ee
for $i=1, \dots, n$ and for every test function $\vp \in W^{1,2}(\Omega)$.

For the discrete play operator \eqref{de4}, we have an explicit expression
\be{de4b}
\xi^r_i(x) = \min\{u_i(x) + r, \max\{\xi^r_{i-1}(x), u_i(x) - r\}\}.
\ee
Formula \eqref{de4b} gives an analytical justification for Figure~\ref{f2}. In addition, using \eqref{de3}, we see that $G[u]_i$ can be represented by a Nemytskii (superposition) operator
$$
G[u]_i(x) = \tilde G_i(x, u_i(x))
$$
with a bounded function $\tilde G_i$ which is continuous and non-decreasing in the second variable. Hence, by the standard theory of semilinear monotone elliptic problems in Hilbert spaces (see, e.\,g., \cite[Theorem~10.49]{renrog}), Eq.~\eqref{de0v} has a unique solution $u_i \in W^{1,2}(\Omega)$ for every $i = 1, \dots, {n}$. As a special case of \cite[Section~2.3.3]{gris}, we also have $u_i \in W^{2,2}(\Omega)$.


\section{Estimates independent of $\tau$}\label{sec:itau}

The next program is to derive below in Subsections~\ref{sec:ubound}--\ref{sec:time} a series of estimates independent of the time step $\tau$, which will allow us to pass to the limit as $\tau \to 0$. We first prove that the solutions $u_i$ to \eqref{de0v} remain bounded in the sup-norm independently of $i$.


\subsection{Uniform upper bound}\label{sec:ubound}

With the intention to prove a uniform upper bound for $u_i$, let us consider the linear elliptic problem
\be{ell1}
\io \nabla v\cdot \nabla \vp \dd x + \ipo b(x) v \vp \dd s(x) = \io \tilde h \vp \dd x
\ee
for every $\vp \in W^{1,2}(\Omega)$ with a given function $\tilde h \in L^\infty(\Omega)$. It has the following property.

\begin{proposition}\label{pm}
Let $\tilde h \in L^\infty(\Omega)$, $|\tilde h(x)|\le K$ a.\,e. Then there exists $V>0$ such that the solution $v$ to \eqref{ell1} satisfies the estimate $|v(x)| \le V$ for a.\,e. $x \in \Omega$.
\end{proposition}

A proof of this fact can be deduced for example from the general theory in \cite[Chapter~2]{pro-wei}. We show here another classical approach based on Moser iterations which we repeat here for the reader's convenience.

\bpf{Proof of Proposition~\ref{pm}}
We choose an arbitrary $R>0$, and for $v \in W^{1,2}(\Omega)$ and $\kappa \ge 0$ put $\bar v = \max\{-R, \min\{v, R\}\}$, $\Gamma_{\kappa}(v) = v|\bar v|^{\kappa}$. The identity
$$
\nabla \Gamma_\kappa(v) = \left\{
\begin{array}{ll}
(\kappa+1)|v|^\kappa \nabla v  & \for |v| < R,\\[2mm]
R^\kappa \nabla v & \for |v| \ge R,
\end{array}
\right.
$$
holds for all $\kappa \ge 0$ and a.\,e. $x\in \Omega$, and we use it for testing \eqref{ell1} by $\vp = \Gamma_{2j}(v)$ with some $j \ge 0$ which will be specified later on. This is indeed an admissible test function, and using the inequality
$$
\nabla v\cdot \nabla \vp = \nabla v\cdot \nabla\Gamma_{2j}(v) \ge \frac{2j+1}{(j+1)^2} \big|\nabla\Gamma_j(v)\big|^2\ge \frac{1}{j+1} \big|\nabla(v |\bar v|^j)\big|^2,
$$
which holds for a.\,e. $x \in \Omega$, we obtain that
\be{ell2}
\frac{1}{j{+}1} \io\big|\nabla(v |\bar v|^j)\big|^2\dd x + \ipo b(x)\big|v |\bar v|^j\big|^2\dd s(x) \le K\io \big|v |\bar v|^{2j}\big| \dd x \le K \io \max\left\{1,\big|v |\bar v|^{j}\big|^2\right\} \dd x.
\ee
We now use the embedding formula $W^{1,2}(\Omega) \hookrightarrow L^p(\Omega)$ for $\frac12 > \frac1p \ge \frac12 - \frac1N$ with embedding constant $C>0$ and obtain that
\be{ell3}
\big|v |\bar v|^{j}\big|_p^2 \le CK(j+1) \max\left\{1,\big|v |\bar v|^{j}\big|_2^2\right\},
\ee
where the symbol $|\cdot|_p$ denotes the norm in $L^p(\Omega)$. Passing to the limit as $R\to \infty$ we deduce from \eqref{ell3} the implication
\be{ell4}
\io|v|^{2(j+1)} \dd x < \infty \  \Longrightarrow \ \io|v|^{p(j+1)} \dd x < \infty.
\ee
We know that $\io|v|^2 \dd x < \infty$. Hence, $\io|v|^q \dd x < \infty$ for every $q\ge 1$, and \eqref{ell3} yields for all $j \ge 0$ that
\be{ell5}
\left(\io\max\{1,|v|\}^{p(j+1)}\dd x\right)^{2/p} \le CK(j+1) \io \max\{1,|v|\}^{2(j+1)}\dd x
\ee
with constants $C>0$, $K>0$ independent of $j$.

For $k \in \nat \cup \{0\}$ we now put
$$
{d} = \frac{p}{2}, \quad j = {d}^k - 1.
$$
Then \eqref{ell5} reads
\be{ell6}
\left(\io\max\{1,|v|\}^{2{d}^{k+1}}\dd x\right)^{1/{d}} \le CK{d}^k \io\max\{1,|v|\}^{2{d}^k}\dd x.
\ee
Put
$$
V_k := \big|\max\{1,|v|\}\big|_{2{d}^{k}}.
$$
Then \eqref{ell6} is of the form
$$
V_{k+1} \le \left(CK {d}^k\right)^{{d}^{-k}/2} V_k,
$$
that is,
\be{ell7}
\log V_{k+1} - \log V_k \le \frac{{d}^{-k}}{2}\left(\log(CK) + k\log{d}\right).
\ee
The sum over $k=0,1,\dots$ of the terms on the right-hand side of \eqref{ell7} as well as $V_0$ are finite. We thus have $\sup_{k\in \nat\cup\{0\}} V_k < \infty$, which proves the statement of Proposition~\ref{pm}.
\epf

Using this result and a discrete variant of the \textit{Hilpert inequality} (see \cite{hilp}, \cite[Section~IX.2]{vis}) derived below in \eqref{dh2}, we now prove the following uniform estimate for the solutions $u_i$ to \eqref{de0v}.

\begin{proposition}\label{p2}
Let Hypotheses~\ref{hyd} and \ref{hym} hold, and let $\{u_0, u_1, \dots, u_n\} \subset W^{2,2}(\Omega)$ satisfy \eqref{de0v}. Then there exists a constant $\bar U>0$ independent of $n$ such that $|u_i(x)| \le \bar U$ a.\,e. for all $i=0,1,\dots, {n}$.
\end{proposition}

\bpf{Proof}
We choose a constant $U^* > |u^*|_\infty$ and put $v^* = v + U^*$, where $v$ is the solution of \eqref{ell1} with some $\tilde h > |h|_\infty$. Taking a bigger $U^*$ if necessary, we may assume that
\be{inil}
v^*(x) > \Lambda \ \ale
\ee
with $\Lambda$ from \eqref{e6b}. Let $G[v^*]_i$ be the output of the Preisach operator \eqref{de3} associated with the constant sequence $v^*_i = v^*$ and with initial memory curve
\be{inivst}
\lambda^*(x,r) = (v^*(x) - r)^+.
\ee
Then $G[v^*]_i = G[v^*]_{i-1}$ for all $i = 1, \dots, {n}$, and we can rewrite \eqref{ell1} in the form
\be{dell1}
\io \left(\frac1\tau(G[v^*]_i - G[v^*]_{i-1})\vp + \nabla v^*\cdot \nabla \vp\right) \dd x + \ipo b(x) (v^*- U^*) \vp \dd s(x) = \io \tilde h \vp \dd x.
\ee
We now subtract \eqref{dell1} from \eqref{de0v} and test the difference by a regularization $\vp = H_\ve(u_i - v^*)$ of the Heaviside function
\be{hev}
H_\ve(z) = \left\{
\begin{array}{ll}
&0 \ \for \ z \le 0,\\[2mm]
& \displaystyle{\frac{z}{\ve}} \ \for \ z \in (0,\ve),\\[3mm]
& 1 \ \for \ z \ge \ve.
\end{array}
\right.
\ee
We have $\nabla(u_i - v^*)\cdot \nabla H_\ve(u_i - v^*) \ge 0$ and $(h_i - \tilde h)H_\ve(u_i - v^*) \le 0$ a.\,e. in $\Omega$, $(u_i - u_i^* - v^* + U^*)H_\ve(u_i - v^*) \ge 0$ a.\,e. on $\partial\Omega$, so that
\be{dh1}
\io \big((G[u]_i - G[u]_{i-1})-(G[v^*]_i - G[v^*]_{i-1})\big)H_\ve(u_i - v^*) \dd x \le 0,
\ee
for every $\ve > 0$. Letting $\ve \to 0$ we get
\be{dh2}
\io \big((G[u]_i - G[u]_{i-1})-(G[v^*]_i - G[v^*]_{i-1})\big)H(u_i - v^*) \dd x \le 0,
\ee
where $H$ is the Heaviside function
\be{heavi}
H(z) = \left\{
\begin{array}{ll}
0 &\ \for z \le 0,\\
1 &\ \for z>0.
\end{array}\right.
\ee
We now proceed by induction. Assume that
\be{vst}
u_{i-1}(x) \le v^*(x), \ G[u]_{i-1}(x) \le G[v^*]_{i-1}(x),\ \xi^r_{i-1}(x) = 0 \ \for r \ge v^*(x) \ \ale
\ee
This is true for $i=1$. Indeed, by \eqref{de5}, \eqref{e6}--\eqref{e6b} and \eqref{inil}--\eqref{inivst} we have $v^*(x) \ge |u_0|_\infty$, $\lambda^*(x,r) \ge \lambda(x,r)$ a.\,e., and \eqref{vst} for $i=1$ follows. If \eqref{vst} holds for some $i>1$, then, by \eqref{dh2},
\be{dh3}
\io \big(G[u]_i -G[v^*]_i\big)H(u_i - v^*) \dd x \le \io \big(G[u]_{i-1} -G[v^*]_{i-1}\big)H(u_i - v^*) \dd x \le 0.
\ee
Assume that there exists a set $\Omega^\sharp \subset \Omega$ of positive measure such that $u_i(x) > v^*(x)$ for $x \in \Omega^\sharp$. Then for $x \in \Omega^\sharp$ we have $\xi_i^r(x) = (u_i(x) - r)^+$ for $r \ge 0$, and
\begin{align*}
G[u]_i(x) &= \int_0^{u_i(x)}\int_0^{u_i(x)-r}\rho(x,r,v)\dd v\dd r \ge \int_0^{v^*(x)}\int_0^{u_i(x) - r}\rho(x,r,v)\dd v\dd r,\\
G[v^*]_i(x) &= \int_0^{v^*(x)}\int_0^{v^*(x) - r}\rho(x,r,v)\dd v\dd r,
\end{align*}
hence,
$$
G[u]_i(x) - G[v^*]_i(x) \ge \int_0^{v^*(x)}\int_{v^*(x)-r}^{u_i(x) - r} \rho(x,r,v)\dd v\dd r  = \int_{v^*(x)}^{u_i(x)} \int_0^{v^*(x)} \rho(x,r,w-r)\dd r\dd w.
$$
Choose in \eqref{e3a} the value $U = u_i(x)$. We get
$$
G[u]_i(x) - G[v^*]_i(x) \ge \rho_0(u_i(x))\Lambda(u_i(x) - v^*(x)) = C_i(x) (u_i(x) - v^*(x))
$$
for $x \in \Omega^\sharp$, with $C_i(x) > 0$. This, together with \eqref{dh3}, contradicts the assumption that $\meas (\Omega^\sharp) > 0$. Hence, $u_i(x) \le v^*(x)$ a.\,e.
From \eqref{de4b} it follows that $\xi^r_i(x) \le \lambda^*(x,r)$ a.\,e. with $\lambda^*$ from \eqref{inivst}, hence $G[u]_i(x) \le G[v^*]_i(x)$ a.\,e., and the induction step is complete. By symmetry we find a constant $C_*>0$ such that $u_i(x) \ge -C_*$ a.\,e., which concludes the proof. 
\epf


\subsection{Initial time estimate}\label{sec:init}

As pointed out in the Introduction, the first issue in letting the time step tend to $0$ is an estimate for the initial time increment independent of the discretization, which does not immediately follow from \eqref{de0v}. To overcome the problem, we use Hypothesis~\ref{hym} to construct a fictitious history of the process $u_{-1}(x)$ with associated memory curve $\lambda_{-1}(x,r) = \xi^r_{-1}(x)$ which is compatible with \eqref{e0e}--\eqref{de1}. Since $u_0 \in W^{2,2}(\Omega)$ by Hypothesis~\ref{hyd}, we only need
\be{e7}
\frac1\tau(G[u]_0 - G[u]_{-1})(x) = \Delta u_0(x) + h(x,0)
\ee
to hold a.\,e. in $\Omega$ together with boundary condition \eqref{c2a}. Then $u_0$ is also a weak solution to \eqref{e7} in the sense of \eqref{de0v} for $i=0$. This will be exploited later in Subsection~\ref{sec:time}. Note that \eqref{e7} also allows us to apply Proposition~\ref{p2} to $u_0$, and obtain
\be{u0}
\supess_{x\in \Omega} |u_0(x)|\le \bar U.
\ee
Compatibility with \eqref{de3} and \eqref{de4b} is guaranteed provided
\begin{align}\label{mi0}
G[u]_{-1}(x) &= \bar G + \int_0^{\infty}\int_0^{\xi^r_{-1}(x)} \rho(x,r,v)\dd v\dd r,\\ \label{mi1}
\xi^r_0(x) &= \min\{u_0(x) + r, \max\{\xi^r_{-1}(x), u_0(x) - r\}\}.
\end{align}

\begin{proposition}\label{min1}
Let Hypothesis~\ref{hym} hold, and let $\bar U>0$ be as in Proposition~\ref{p2} and such that $r_0(x) \le \bar U$ a.\,e. With the notation of \eqref{e3a} put $\rho^* = \rho_0(3\bar U)$. Then for
\be{tz}
\tau < \tau_0 := \frac{\rho^*}{2L^2},
\ee
there exists $u_{-1} \in L^\infty(\Omega)$ and a memory curve $\lambda_{-1}$ satisfying the conditions of Definition~\ref{dpr} and such that the identities \eqref{e7}, \eqref{mi0}--\eqref{mi1} are satisfied a.\,e. in $\Omega$ and all $r>0$ with $\lambda_{-1}(x,r) = \xi^r_{-1}(x)$, $\lambda_{-1}(x,0) = u_{-1}(x)$ a.\,e., and
\be{umin}
\frac1\tau|u_0(x) - u_{-1}(x)| \le C \ \ \ale
\ee
with a constant $C>0$ independent of $\tau$ and $x$.
\end{proposition}

\begin{figure}[htb]
\begin{center}
\includegraphics[width=14cm]{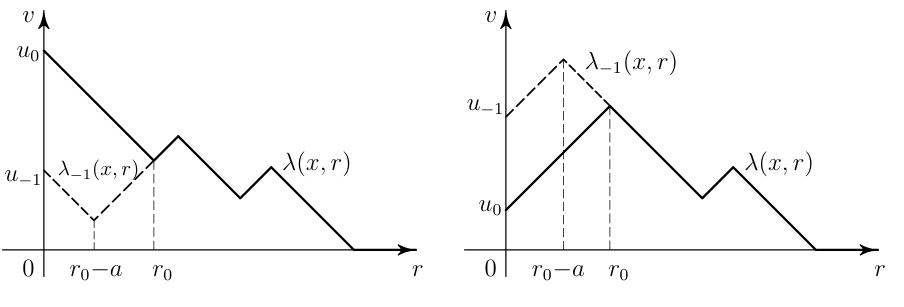}
\caption{Backward time step in $\Omega_+$ (left) and in $\Omega_-$ (right).} \label{f3}
\end{center}
\end{figure}

\bpf{Proof}
We define the sets
$$
\begin{array}{ll}
\Omega_+ &= \{x \in \Omega: \Delta u_0(x) + h(x,0) >0\},\\
\Omega_- &= \{x \in \Omega: \Delta u_0(x) + h(x,0) <0\},\\
\Omega_0 &= \{x \in \Omega: \Delta u_0(x) + h(x,0) =0\}.
\end{array}
$$
Consider first $x \in \Omega_+$. We have by \eqref{c2} and \eqref{de5} that
$$
\lambda(x,r) = u_0(x) - r \ \for r\in [0,r_0(x)].
$$
We choose $a(x) \in [0,r_0(x)/2]$, define $\lambda_{-1}(x,r)$ by the formula
$$
\lambda_{-1}(x,r) = \left\{
\begin{array}{ll}
u_0(x) - r - 2{a(x)} &\ \for\ r\in [0, r_0(x) - {a(x)}],\\
u_0(x) - 2r_0(x) + r&\ \for\ r\in (r_0(x) - {a(x)}, r_0(x)),\\
\lambda(x,r) &  \ \for\ r \ge r_0(x).
\end{array}
\right.
$$
It is easy to see that $\lambda_{-1}$ satisfies the conditions of Definition~\ref{dpr}. Put for $x \in \Omega_+$
$$
u_{-1}(x) = \lambda_{-1}(x,0), \quad G[u]_{-1}(x) = \bar{G} + \int_0^\infty\int_0^{\lambda_{-1}(x,r)}\rho(x,r,v)\dd v\dd r, \quad \xi^r_{-1}(x) = \lambda_{-1}(x,r)
$$
in agreement with \eqref{mi0}, see Figure~\ref{f3}. Moreover,
$$
\xi^r_{-1}(x) = \lambda_{-1}(x,r) \le \lambda(x,r) = u_0(x) - r \quad \ale \ \mbox{for all } r \in [0,r_0(x)],
$$	
so that, recalling \eqref{de5}, condition \eqref{mi1} is satisfied.
We have $2{a(x)} = u_0(x) - u_{-1}(x)$ and
\begin{align}\label{one}
G[u]_0(x) - G[u]_{-1}(x) &= \int_0^\infty\int_{\lambda_{-1}(x,r)}^{\lambda(x,r)}\rho(x,r,v)\dd v\dd r\\ \nonumber
&= \int_0^{r_0(x) - {a(x)}}\int_{u_0(x) - r - 2{a(x)}}^{u_0(x) - r}\rho(x,r,v)\dd v\dd r\\ \nonumber
&\quad + \int_{r_0(x) - {a(x)}}^{r_0(x)}\int_{u_0(x) - 2r_0(x) + r}^{u_0(x) - r}\rho(x,r,v)\dd v\dd r.
\end{align}

In agreement with $\eqref{GB}$, we denote the right-hand side of \eqref{one} by $B(x,{a(x)})$, where
$$
B(x,{a}) = \int_0^{r_0(x) - {a}}\int_{u_0(x) - r - 2{a}}^{u_0(x) - r}\rho(x,r,v)\dd v\dd r
+ \int_{r_0(x) - {a}}^{r_0(x)}\int_{u_0(x) - 2r_0(x) + r}^{u_0(x) - r}\rho(x,r,v)\dd v\dd r.
$$
By assumption \eqref{e3a} and estimate \eqref{u0}, and recalling that $r_0(x) \le \bar U$ a.\,e. and $\rho^* = \rho_0(3\bar U)$, we have 
$$
B_a(x,a) := \frac{\partial}{\partial{a}} B(x,{a}) = 2\int_0^{r_0(x) - {a}}\rho(x,r,u_0(x) - r -2{a})\dd r
\ge 2\rho^*(r_0(x) - {a}).
$$
We have indeed $B(x,0) = 0$, and $B_a(x,{a}) \ge \rho^* r_0(x)$ for ${a} \le r_0(x)/2$. Hence,
\be{gd}
B(x,{a}) \ge \rho^* r_0(x){a} \ \for \ 0 < {a} \le \frac{r_0(x)}{2}.
\ee
In particular, for
\be{hd}
\hat a(x) := \frac{r_0(x)}{2}
\ee
 we have
\be{delta}
B(x,\hat a(x)) \ge \frac{\rho^*}{2} r_0^2(x)\ge  \frac{\rho^*}{2L^2} (\Delta u_0(x) + h(x,0)).
\ee 
Again, since $B(x,0) = 0$ and $a \mapsto B(x,a)$ is an increasing function, we conclude from \eqref{one} and \eqref{delta} that for each $\tau>0$ satisfying \eqref{tz} there exists ${a(x)} \in (0,\hat a(x))$ such that
\be{two}
G[u]_0(x) - G[u]_{-1}(x) = B(x,{a(x)}) = \tau(\Delta u_0(x) + h(x,0)),
\ee
which is precisely \eqref{e7}. Finally, by \eqref{c1}, \eqref{gd}, and \eqref{two}, we have
$$
u_0(x) - u_{-1}(x) = 2{a(x)} \le \frac{2B(x,{a(x)})}{\rho^* r_0(x)} \le \frac{2L\tau }{\rho^*} \sqrt{\Delta u_0(x) + h(x,0)} \le 
\frac{2L^2\Lambda\,\tau}{\rho^*},
$$
so that \eqref{umin} holds with $C = 2L^2\Lambda/\rho^*$.

For $x \in \Omega_-$ we proceed similarly. We have by \eqref{c2} and \eqref{de5} that
$$
\lambda(x,r) = u_0(x) + r \ \for r\in [0,r_0(x)],
$$
and putting
$$
\lambda_{-1}(x,r) = \left\{
\begin{array}{ll}
u_0(x) + r + 2{a(x)} &\ \for\ r\in [0, r_0(x) - {a(x)}],\\
u_0(x) + 2r_0(x) - r&\ \for\ r\in (r_0(x) - {a(x)}, r_0(x)),\\
\lambda(x,r) &  \ \for\ r \ge r_0(x),
\end{array}
\right.
$$
see Figure~\ref{f3}, we obtain \eqref{e7} and \eqref{umin} similarly as above. For $x\in \Omega_0$ we simply put $\lambda_{-1}(x,r) = \lambda(x,r)$ and $u_{-1}(x) = u_0(x)$, and the assertion holds trivially.
\epf


\subsection{Estimates for the space derivatives}\label{sec:espace}

Let us first notice that by hypothesis we have
$$
h \in L^\infty(\Omega\times (0,T)) \cap W^{1,2}(0,T; L^2(\Omega)), \quad u^* \in L^\infty(\partial\Omega\times (0,T)) \cap W^{1,2}(\partial\Omega\times (0,T)),
$$ 
hence,
\be{hi0}
\max_{i=0,\dots,{n}}\supess_{x\in \Omega} |h_i(x)|\le C, \quad \max_{i=0,\dots,{n}}\supess_{x\in \partial\Omega} |u^*_i(x)|\le C,
\ee
and
\be{hi}
\frac1\tau \sum_{i=1}^{n}\io |h_i(x) - h_{i-1}(x)|^2\dd x \le C, \quad \frac1\tau \sum_{i=1}^{n}\ipo |u^*_i(x) - u^*_{i-1}(x)|^2\dd s(x) \le C
\ee
with a constant $C>0$ independent of $\tau$. Moreover, from the validity of \eqref{de0v} for $i = 0$ stated in \eqref{e7}, we can apply Proposition~\ref{p2} to $\{u_0,u_1,\dots,u_n\}$ and obtain
\be{ui}
\max_{i=0,\dots,{n}}\supess_{x\in \Omega} |u_i(x)|\le \bar U.
\ee

Testing \eqref{de0v} by $\vp = u_i$ and using \eqref{ene1} we get
$$
\frac1\tau \io(E[u]_i - E[u]_{i-1})\dd x + \io |\nabla u_i|^2\dd x + \ipo b(x) (u_i - u^*_i)u_i\dd s(x) \le \io h_i u_i \dd x,
$$ 
hence, summing up over $i=1, \dots, {n}$ and invoking \eqref{iniene}, \eqref{hi0}, and \eqref{ui}
\be{es0}
\tau\sumj \left(\io |\nabla u_i|^2\dd x + \ipo b(x) |u_i|^2\dd s(x)\right) \le C
\ee
with a constant $C>0$ independent of $\tau$.

The next standard step consists in testing \eqref{de0v} by $u_i - u_{i-1}$ as a time discrete counterpart of $u_t$, which yields
\begin{align*}
&\frac1\tau \io(G[u]_i - G[u]_{i-1})(u_i - u_{i-1})\dd x + \io \nabla u_i\cdot \nabla(u_i - u_{i-1})\dd x + \ipo b(x) u_i(u_i - u_{i-1})\dd s(x)\\
&\hspace{10mm} = \io h_i (u_i - u_{i-1}) \dd x+\ipo b(x) u^*_i (u_i - u_{i-1})\dd s(x).
\end{align*} 
On the left-hand side of this identity we have
$$
\nabla u_i\cdot \nabla(u_i - u_{i-1}) \ge \frac12\left(|\nabla u_i|^2 - |\nabla u_{i-1}|^2\right), \quad u_i(u_i - u_{i-1}) \ge \frac12\left(|u_i|^2 - |u_{i-1}|^2\right),
$$
and we rewrite the terms on the right-hand side as
$$
h_i(u_i - u_{i-1}) = h_i u_i - h_{i-1} u_{i-1} + (h_i - h_{i-1})u_{i-1}, 
$$
and similarly
$$
u^*_i(u_i - u_{i-1}) = u^*_i u_i - u^*_{i-1} u_{i-1} + (u^*_i - u^*_{i-1})u_{i-1}.
$$
Put
$$
V[u]_i = \io \left(\frac12 |\nabla u_i|^2 - h_i u_i\right)\dd x + \ipo b(x)\left(\frac12 |u_i|^2 - u^*_i u_i \right)\dd s(x).
$$
The above computations then yield that
\begin{align}\nonumber
&\hspace{-8mm}\frac1\tau \io(G[u]_i - G[u]_{i-1})(u_i - u_{i-1})\dd x + V[u]_i - V[u]_{i-1}\\ \label{aux1}
& \le \io (h_i - h_{i-1})u_{i-1}\dd x + \ipo b(x)(u^*_i - u^*_{i-1})u_{i-1}\dd s(x).
\end{align}
From Hypothesis~\ref{hyd} for $u^0=u_0$ and estimates \eqref{hi0}--\eqref{ui} we thus get
\be{es1}
\frac1\tau \sumj\io(G[u]_i - G[u]_{i-1})(u_i - u_{i-1})\dd x + \max_{i=1, \dots, {n}} \left(\io |\nabla u_i|^2\dd x + \ipo b(x) |u_i|^2 \dd s(x)\right) \le C
\ee
with a constant $C>0$ independent of $\tau$. Finally, by comparison in \eqref{e0e}, inequality \eqref{de5a} yields that
\be{es2}
\tau\sumj\io |\Delta u_i|^2\dd x \le C.
\ee
The regularity of $\partial\Omega$ and of the boundary data guarantee that the results of \cite[Section~2.3.3]{gris} can be applied here to check that
\be{es2b}
|u|_{2,2,b} := \left(\io |\Delta u|^2 \dd x + \ipo b(x) |u|^2 \dd s(x)\right)^{1/2}
\ee
is an equivalent norm in $W^{2,2}(\Omega)$, and by \eqref{es0} and \eqref{es2} we have
\be{es2a}
\tau\sumj|u_i|_{2,2,b}^2 \le C
\ee
with a constant $C>0$ independent of $\tau$.

The estimates \eqref{es0}, \eqref{es1}, \eqref{es2a} are indeed not sufficient for passing to the limit as ${n} \to \infty$, as we do not control the time increments $u_i - u_{i-1}$ of $u$. We have to use a tool which is specific for systems with hysteresis, namely a convexity argument. This is where the convexifiability of $G$ in Theorem~\ref{t1} comes into play.


\subsection{Convexity}\label{sec:conr}

Consider a uniformly counterclockwise convex Preisach operator $P$ on an interval $[-U,U]$ in the sense of Definition~\ref{dpc}, in discrete form similarly as in \eqref{de3}, that is,
\be{pe3}
P[w]_i(x,t) = \bar{G} +  \int_0^{\infty}\int_0^{\zeta^r_i(x)} \sigma(x,r,v)\dd v\dd r.
\ee
Here, as in \eqref{de4}--\eqref{de5}, $\zeta^r_i = \play_r[\mu,w]_i$ is the output of the discrete play with input $\{w_i: i=-1,0,\dots n\} \subset L^\infty(\Omega)$ and initial memory configuration $\mu(x,r)$ satisfying
\be{emm}
|\mu(x,r_2) - \mu(x, r_1)| \le |r_2 - r_1|, \quad \mu(x,r) = 0 \ \for r\ge U\ \ale
\ee
Let us assume furthermore that $\{w_i\}$ is such that $|w_i(x)| \le U$ a.\,e. for all $i=-1,0,\dots n$, and let $i\ge 0$ be fixed. Since $P$ is by hypothesis uniformly counterclockwise convex on $[-U,U]$, there exist increasing functions $B_+^{i-1}(x,\cdot): [w_{i-1}(x), U] \to \real$, $B_-^{i-1}(x,\cdot): [-U, w_{i-1}(x)]\to \real$, $B_+^{i-1}(x,\cdot)$ uniformly convex, $B_-^{i-1}(x,\cdot)$ uniformly concave of the form
\begin{align}\label{b+}
B_+^{i-1}(x,w) &= \int_0^{r^{i-1}(x,w)}\int_{\zeta^r_{i-1}(x)}^{w-r} \sigma(x,r,v)\dd v\dd r \ \ \for\ w \in [w_{i-1}(x), U],\\ \label{b-}
B_-^{i-1}(x,w) &= -\int_0^{r_{i-1}(x,w)}\int^{\zeta^r_{i-1}(x)}_{w+r} \sigma(x,r,v)\dd v\dd r \ \ \for\ w \in [-U, w_{i-1}(x)],
\end{align}
with $r^{i-1}(x,w) = \min\{r>0: w-r \le \zeta^r_{i-1}(x)\}$, $r_{i-1}(x,w) = \min\{r>0: w+r \ge \zeta^r_{i-1}(x)\}$, and by \eqref{GB} we have
\begin{align}\label{p+}
P[w]_i(x) &= P[w]_{i-1}(x) + B_+^{i-1}(x,w_i) \ \ \for \ w_i(x) \ge w_{i-1}(x),\\ \label{p-}
P[w]_i(x) &= P[w]_{i-1}(x) + B_-^{i-1}(x,w_i) \ \ \for \ w_i(x) < w_{i-1}(x). 
\end{align}
The functions $B_+^{i-1}, B_-^{i-1}$ are the ascending and descending branches of $P$ at discrete time $i-1$, respectively.
For simplicity, we do not explicitly display in the next developments the dependence of $B_{\pm}^{i-1}$ on $x$ and $i$, and we write simply $B_\pm(w)$ instead of $B_\pm^{i-1}(x,w)$.

The hypothesis of uniform counterclockwise convexity of $P$  on $[-U,U]$ means that there exists $\beta>0$ such that all ascending branches $B_+: [\hat a, \hat b] \subset [-U, U] \to \real$ are uniformly convex in the sense
\be{w+}
B'_+(w) >0, \quad \frac{B'_+(w_1) - B'_+(w_2)}{w_1 - w_2} \ge \beta \  \ \mbox{ for a.\,e. } w, w_1\ne w_2 \in (\hat a, \hat b),
\ee
and all descending branches $B_-: [\hat a, \hat b] \subset [-U, U] \to \real$ are uniformly concave in the sense
\be{w-}
B'_-(w) > 0, \quad \frac{B'_-(w_1) - B'_-(w_2)}{w_1 - w_2} \le -\beta  \ \ \mbox{ for a.\,e. } w, w_1 \ne w_2 \in (\hat a, \hat b).
\ee

\begin{lemma}\label{cl1}
Let $B_+, B_-$ satisfy \eqref{w+}--\eqref{w-}. Then for every $-U \le \hat a \le a \le b \le c \le \hat b \le U$ we have
$$
\frac{B_+(c) - B_+(b)}{c-b} - \frac{B_+(b) - B_+(a)}{b-a} \ge \frac{\beta}{2}(c-a),
$$
and
$$
\frac{B_-(c) - B_-(b)}{c-b} - \frac{B_-(b) - B_-(a)}{b-a} \le -\frac{\beta}{2}(c-a),
$$
where for $c=b$ the quotients $(B_+(c) - B_+(b))/(c-b)$ and $(B_-(c) - B_-(b))(c-b)$ are to be interpreted as $B'_+(c-)$ and $B'_-(c-)$, and for $b=a$ the quotients $(B_+(b) - B_+(a))/(b-a)$ and $(B_-(b) - B_-(a))(b-a)$ are to be interpreted as $B'_+(a+)$ and $B'_-(a+)$.
\end{lemma}

\bpf{Proof}
We prove only the statement for $B_+$, the other one is obtained by symmetry. We have
$$
\frac{B_+(c) - B_+(b)}{c-b} - \frac{B_+(b) - B_+(a)}{b-a} = \frac1{c-b}\int_b^c B'_+(v) \dd v - \frac1{b-a} \int_a^b B'_+(w) \dd w.
$$
By substitution $v = (c-b)z + b$ and $w = (b-a)z + a$ for $z \in (0,1)$, we rewrite the integrals on the right-hand side as
$$
\int_0^1 \big(B'_+((c-b)z + b) - B'_+((b-a)z + a)\big)\dd z \ge \beta \int_0^1 \big((c+a - 2b)z + (b-a)\big)\dd z = \frac{\beta}{2}(c-a),
$$
which we wanted to prove.
\epf

\begin{lemma}\label{cl2}
Let $B_+, B_-$  be as in Lemma~\ref{cl1}, let $f:\real \to \real$ be an increasing function such that $f(0) = 0$, and put $F(w) = \int_0^w f(v)\dd v$. Then for all $-U \le \hat a\le a\le b\le c\le \hat b \le U$  we have
\begin{align*}
(B_+(c) - 2B_+(b) + B_+(a))f(c-b) &\ge \frac{B_+(c) - B_+(b)}{c-b} F(c-b) - \frac{B_+(b) - B_+(a)}{b-a} F(b-a)\\[2mm]
&\quad + \frac{\beta}2(c-a)\big((c-b)f(c-b) - F(c-b)\big),
\end{align*}
and for all $U \ge \hat b\ge a\ge b\ge c\ge \hat a \ge -U$
\begin{align*}
(B_-(c) - 2B_-(b) + B_-(a))f(c-b) &\ge \frac{B_-(c) - B_-(b)}{c-b} F(c-b) - \frac{B_-(b) - B_-(a)}{b-a} F(b-a)\\[2mm]
&\quad + \frac{\beta}2(a-c)\big((c-b)f(c-b) - F(c-b)\big).
\end{align*}
\end{lemma}

\bpf{Proof}
We prove again only the statement for $B_+$, the other one is obtained by replacing $B_-(v)$ with $-B_+(-v)$ and $f(v)$ with $- f(-v)$. We have
\begin{align*}
&(B_+(c) - 2B_+(b) + B_+(a))f(c-b) - \frac{B_+(c) - B_+(b)}{c-b} F(c-b) + \frac{B_+(b) - B_+(a)}{b-a} F(b-a)\\[2mm]
&\quad =\left(\frac{B_+(c) - B_+(b)}{c-b} - \frac{B_+(b) - B_+(a)}{b-a}\right)\big((c-b)f(c-b) - F(c-b)\big)\\
&\qquad + \frac{B_+(b) - B_+(a)}{b-a}\big(f(c-b)\big((c-b) - (b-a)\big) - \big(F(c-b) - F(b-a)\big)\big).
\end{align*}
Since $F$ is convex and $B_+$ is increasing, we have
$$
\frac{B_+(b) - B_+(a)}{b-a}\big(f(c-b)\big((c-b) - (b-a)\big) - \big(F(c-b) - F(b-a)\big)\big) \ge 0,
$$
and the assertion follows from Lemma~\ref{cl1}.
\epf

The following result is crucial for proving the existence of solutions to Problem \eqref{e0}.

\begin{proposition}\label{pc}
Let $\{w_i(x): i=-1,0,\dots, {n}\}$ be a sequence in $L^\infty(\Omega)$, $|w_i(x)| \le U$ a.\,e. for all $i = -1,0,\dots, {n}$. Let $P$ in \eqref{pe3} be uniformly counterclockwise convex on $[-U,U]$, with ascending branches $B_+$ and descending branches $B_-$ satisfying \eqref{w+} and \eqref{w-}, respectively. Let $f, F$ be as in Lemma~\ref{cl2}. Then for a.\,e. $x \in \Omega$ we have
\be{Pi}
\begin{aligned}
&\sumi (P[w]_{i+1} - 2P[w]_i + P[w]_{i-1})f(w_{i+1} - w_i) + \frac{P[w]_0 - P[w]_{-1}}{w_0 - w_{-1}}F(w_0 - w_{-1}) \\[2mm]
&\qquad \ge \frac{\beta}{2}\sumi |w_{i+1} - w_i|\big((w_{i+1} - w_i) f(w_{i+1} - w_i) - F(w_{i+1} - w_i)\big).
\end{aligned}
\ee
\end{proposition}

\bpf{Proof}
We choose an arbitrary $i \in \{0,1,\dots, {n-1}\}$ and distinguish the cases
\begin{itemize}
\item[(i)] $w_{i-1} \le w_i \le w_{i+1}$

\item[(ii)] $w_{i-1} \ge w_i \ge w_{i+1}$

\item[(iii)] $w_{i-1} > w_i, w_{i} < w_{i+1}$

\item[(iv)] $w_{i-1} < w_i, w_i > w_{i+1}$
\end{itemize}

In case (i) we first check that the values of $w_{i-1}\le w_i \le w_{i+1}$ are situated on the same ascending branch $B_+ = B_+^{i-1}$. We present here a formal proof of this fact which can be visualized on Figure~\ref{f2}. By \eqref{p+} we obtain
\begin{align*}
P[w]_{i+1}(x) &= P[w]_{i}(x) + B_+^{i}(w_{i+1}) \\
&= P[w]_{i-1}(x) + B_+^{i-1}(w_i) + B_+^{i}(w_{i+1}),
\end{align*}
so that, by \eqref{GB}, we only have to prove that for $w>w_i$ we have
\be{i+1}
B_+^{i-1}(w_i) + B_+^{i}(w) = B_+^{i-1}(w).
\ee
This will be achieved by exploiting the definition \eqref{b+}. First of all, note that $r^{i-1}(w) = r^{i}(w)$. Indeed, let us introduce the sets
\begin{align*}
	J_{i}(w) &= \{r > 0: w-r \le \zeta_i^r\},\\
	J_{i-1}(w) &= \{r > 0: w-r \le \zeta_{i-1}^r\},\\
	J_{i-1}(w_i) &= \{r > 0: w_i-r \le \zeta_{i-1}^r\}.
\end{align*}
From the inequality $w > w_i$, it follows that
\begin{equation}\label{J1}
	J_{i-1}(w) \subset J_{i-1}(w_i).
\end{equation}
Moreover, for $r\in J_{i-1}(w_i)$ we have by \eqref{de4b} that $\zeta_i^r=\zeta_{i-1}^r$. From \eqref{J1} we thus have $\zeta_i^r=\zeta_{i-1}^r$ for all $r\in J_{i-1}(w)$, so that $J_{i-1}(w) \subset J_i(w)$. On the other hand, if there exists $r^* \in J_i(w)\setminus J_{i-1}(w)$, then we would have again by \eqref{de4b} that $w_i - r^* = \zeta_i^{r^*} \ge w-r^*$, which is a contradiction. We thus have $J_{i-1}(w) = J_i(w)$, hence $r^{i-1}(w) = \min J_{i-1}(w) = \min J_{i}(w) = r^{i}(w)$, which is what we wanted to prove. Invoking \eqref{b+} and considering only the integration domains, we thus have
$$
B_+^{i}(w) = \int_0^{r^{i}(w)} \!\!\int_{\zeta_i^r}^{w-r} = \int_0^{r^{i-1}(w_i)}\!\! \int_{w_i-r}^{w-r} + \int_{r^{i-1}(w_i)}^{r^{i-1}(w)}\! \int_{\zeta_{i-1}^r}^{w-r},
$$
hence,
$$
B_+^{i-1}(w_i) + B_+^{i}(w) = \int_0^{r^{i-1}(w_i)}\!\! \int_{\zeta_{i-1}^r}^{w_i-r}+\int_0^{r^{i}(w)} \!\!\int_{\zeta_i^r}^{w-r} = \int_0^{r^{i-1}(w)}\! \int_{\zeta_{i-1}^r}^{w-r} = B_+^{i-1}(w)
$$
and \eqref{i+1} follows, so that we can write
\begin{equation}\label{i+-1}
	P[w]_{i+1} = P[w]_{i-1} + B_+^{i-1}(w_{i+1}).
\end{equation}
Note that \eqref{GB} yields $B_+^{i-1}(w_{i-1}) = 0$, hence, by \eqref{p+} and \eqref{i+-1},
$$
\begin{aligned}
P[w]_{i+1} - P[w]_i &= B_+^{i-1}(w_{i+1}) - B_+^{i-1}(w_i), \\
P[w]_i - P[w]_{i-1} &= B_+^{i-1}(w_i) - B_+^{i-1}(w_{i-1}), \\
P[w]_{i+1} - 2P[w]_i + P[w]_{i-1} &= B_+^{i-1}(w_{i+1}) - 2B_+^{i-1}(w_i) + B_+^{i-1}(w_{i-1}).
\end{aligned}
$$
We can therefore apply Lemma~\ref{cl2} to $B_+^{i-1}$ with $a = w_{i-1}$, $b = w_i$, $c= w_{i+1}$, and obtain
\begin{align}\nonumber
&(P[w]_{i+1} - 2P[w]_i + P[w]_{i-1})f(w_{i+1} - w_i)\\[2mm] \nonumber
&\quad \ge \frac{P[w]_{i+1} - P[w]_i}{w_{i+1}-w_i} F(w_{i+1}-w_i) - \frac{P[w]_i - P[w]_{i-1}}{w_i-w_{i-1}} F(w_i-w_{i-1})\\[2mm] \label{ci}
&\qquad + \frac{\beta}2(w_{i+1}-w_{i-1})\big((w_{i+1}-w_{i})f(w_{i+1}-w_i) - F(w_{i+1}-w_i)\big).
\end{align}
Case (ii) is similar. The values of $w_{i-1} \ge w_i \ge w_{i+1}$ are on the same descending branch $B_-$, and by Lemma~\ref{cl2} we have
\begin{align}\nonumber
&(P[w]_{i+1} - 2P[w]_i + P[w]_{i-1})f(w_{i+1} - w_i)\\[2mm] \nonumber
&\quad \ge \frac{P[w]_{i+1} - P[w]_i}{w_{i+1}-w_i} F(w_{i+1}-w_i) - \frac{P[w]_i - P[w]_{i-1}}{w_i-w_{i-1}} F(w_i-w_{i-1})\\[2mm] \label{cii}
&\qquad + \frac{\beta}2(w_{i-1}-w_{i+1})\big((w_{i+1}-w_{i})f(w_{i+1}-w_i) - F(w_{i+1}-w_i)\big).
\end{align}

In Case (iii), $w_{i-1} > w_i$ are on a descending branch $B_-$, $w_i < w_{i+1}$ are on an ascending branch $B_+$. We have
\begin{align*}
&(P[w]_{i+1} - 2P[w]_i + P[w]_{i-1})f(w_{i+1} - w_i)\\[2mm] 
&\qquad - \frac{P[w]_{i+1} - P[w]_i}{w_{i+1}-w_i} F(w_{i+1}-w_i) + \frac{P[w]_i - P[w]_{i-1}}{w_i-w_{i-1}} F(w_i-w_{i-1})\\[2mm] 
& \quad \ge \frac{P[w]_{i+1} - P[w]_i}{w_{i+1}-w_i}\left((w_{i+1}-w_i)(f(w_{i+1} - w_i)) - F(w_{i+1}-w_i)\right).
\end{align*}
Since $w_i$ is a turning point, we have $B_+'(w_i+) = 0$, and
\begin{align*}
\frac{P[w]_{i+1} - P[w]_i}{w_{i+1}-w_i} &= \frac{B_+(w_{i+1}) - B_+(w_i)}{w_{i+1}-w_i}= \frac{1}{w_{i+1}-w_i}\int_{w_i}^{w_{i+1}} B_+'(z) \dd z\\
& = \frac{1}{w_{i+1}-w_i}\int_{w_i}^{w_{i+1}} (B_+'(z) - B_+'(w_i+))\dd z
\ge \frac\beta{2} (w_{i+1} - w_i),
\end{align*}
where we have used \eqref{w+}, hence
\begin{align}\nonumber
&(P[w]_{i+1} - 2P[w]_i + P[w]_{i-1})f(w_{i+1} - w_i)\\[2mm] \nonumber
&\quad \ge \frac{P[w]_{i+1} - P[w]_i}{w_{i+1}-w_i} F(w_{i+1}-w_i) - \frac{P[w]_i - P[w]_{i-1}}{w_i-w_{i-1}} F(w_i-w_{i-1})\\[2mm] \label{ciii}
&\qquad + \frac\beta{2}(w_{i+1}-w_i)\big((w_{i+1}-w_{i})f(w_{i+1}-w_i) - F(w_{i+1}-w_i)\big).
\end{align}
Finally, in case (iv), $w_{i-1} < w_i$ are on an ascending branch $B_+$, $w_i > w_{i+1}$ are on a descending branch $B_-$, and we have as in case (iii)
\begin{align*}
&(P[w]_{i+1} - 2P[w]_i + P[w]_{i-1})f(w_{i+1} - w_i)\\[2mm] 
&\qquad - \frac{P[w]_{i+1} - P[w]_i}{w_{i+1}-w_i} F(w_{i+1}-w_i) + \frac{P[w]_i - P[w]_{i-1}}{w_i-w_{i-1}} F(w_i-w_{i-1})\\[2mm] 
& \quad \ge \frac{P[w]_{i+1} - P[w]_i}{w_{i+1}-w_i}\left((w_{i+1}-w_i)(f(w_{i+1} - w_i)) - F(w_{i+1}-w_i)\right).
\end{align*}
Since $w_i$ is again a turning point, we have $B_-'(w_i+) = 0$, and
\begin{align*}
\frac{P[w]_{i+1} - P[w]_i}{w_{i+1}-w_i} &= \frac{B_-(w_{i+1}) - B_-(w_i)}{w_{i+1}-w_i}= \frac{1}{w_{i}-w_{i+1}}\int^{w_i}_{w_{i+1}} B_-'(z) \dd z\\
& = \frac{1}{w_i-w_{i+1}}\int^{w_i}_{w_{i+1}} (B_-'(z) - B_-'(w_i+))\dd z
\ge \frac\beta{2} (w_{i} - w_{i+1}),
\end{align*}
where we have used \eqref{w-}, and we can argue as in case (iii).

In all cases, summing up over $i= 0,1,\dots,n-1$ and observing that $F(w_n-w_{n-1})(P[w]_n-P[w]_{n-1})/(w_n-w_{n-1}) \ge 0$, we conclude the proof.
\epf


\subsection{Estimates of the time increment}\label{sec:time}

Assume now that $G$ in \eqref{de3} is convexifiable, that is, for every $U>0$ there exists a function $g:[-U, U] \to [-U, U]$ satisfying \eqref{hg} such that $G$ is of the form
\be{ne0}
G = P \circ g,
\ee
where $P$ is a uniformly counterclockwise convex Preisach operator on $[-U, U]$. Let us choose $U=\bar U$ with $\bar U$ from Proposition~\ref{p2} and the corresponding function $g$. Then \eqref{de0v}, extended to $i=0$ by virtue of \eqref{e7} and \eqref{c2a}, has the form
\be{dp1}
\io \left(\frac1\tau(P[w]_i - P[w]_{i-1})\vp + \nabla u_i\cdot\nabla\vp\right)\dd x + \ipo b(x)(u_i - u^*_i)\vp \dd s(x) = \io h_i\vp\dd x
\ee
with $w_i = g(u_i)$, for $i=0,1,\dots,n$ and for an arbitrary test function $\vp \in W^{1,2}(\Omega)$. The final goal of this section is to exploit the convexity of $P$ to obtain an estimate of the form
$$
\tau^{1-\gamma}\sumi\io |u_{i+1} - u_i|^\gamma\dd x  \le C
$$
for some $\gamma > 1$, with a constant $C>0$ independent of $\tau$. This will be achieved by iterating the exponent $\gamma$, with initial step $\gamma = 1$. The strategy is as explained in Section~\ref{sec:stat}. We take advantage of the convexity using Proposition~\ref{pc}, where we choose successively suitable increasing functions $f_\gamma:\real \to \real$ such that $f_\gamma(0) = 0$, consider the difference of \eqref{dp1} taken at `times' $i+1$ and $i$, and test it by $f_\gamma(w_{i+1}- w_i)$.

The initial step of the iteration scheme is the following.

\begin{proposition}\label{pl1}
The solution to \eqref{dp1} satisfies the estimate
$$
\sumi \io |u_{i+1} - u_i|\dd x \le C
$$
with a constant $C>0$ independent of $\tau$.
\end{proposition}

\bpf{Proof}
Choose
\be{dp2}
f(w) := f_1(w) = \frac{w}{\tau + |w|}
\ee
and put $F(w) = \int_0^w f(v)\dd v = |w| - \tau(\log(\tau + |w|) - \log\tau)$. Then
\be{dp3}
\tilde f(w) := wf(w) - F(w) = \tau\left(\log(\tau+|w|) - \log \tau - \frac{|w|}{\tau+|w|}\right).
\ee
Since $g : [-U,U] \to [-U,U]$ and $w_i = g(u_i)$, we can apply Propositions~\ref{p2} and \ref{pc} with $U = \bar U$, obtaining
\begin{align*}
\frac1\tau \sumi &(P[w]_{i+1} - 2P[w]_i + P[w]_{i-1})f(w_{i+1} - w_i) + \frac{1}{\tau}\frac{P[w]_0-P[w]_{-1}}{w_0-w_{-1}}F(w_0-w_{-1}) \\
&\ge \frac{\beta}{2\tau}\sumi |w_{i+1} - w_i|\tilde f(w_{i+1} - w_i).
\end{align*}
Note that the construction in Subsection~\ref{sec:init} entails the validity of \eqref{Pi} also for $i=0$. Indeed, convexity is preserved as long as the next time step is obtained from the previous one by the play update.

Now, for $|w| \ge \tau(\expe^2 - 1)$ we have $\tilde f(w) \ge \tau$, and for $|w| < \tau(\expe^2 - 1)$ we have $\tilde f(w) \ge 0$. We denote by $J$ the set of all $i \in \{0,1,\dots, n-1\}$ such that $|w_{i+1} - w_i| \ge \tau(\expe^2 - 1)$, and by $J^\perp$ its complement. We therefore have
\begin{align}\label{dp3a}
\frac1\tau\sum_{i\in J} |w_{i+1} - w_i|\tilde f(w_{i+1} - w_i) &\ge \sum_{i\in J} |w_{i+1} - w_i|, \\ \label{dp3c}
\frac1\tau\sum_{i\in J^\perp} |w_{i+1} - w_i|\tilde f(w_{i+1} - w_i) &\ge 0 > \sum_{i\in J^\perp} \big(|w_{i+1} - w_i| - \tau(\expe^2 - 1)\big).
\end{align}
Moreover, for $w_0 \neq w_{-1}$ we have
\be{dp3b}
0 < \frac{F(w_0-w_{-1})}{|w_0-w_{-1}|} \le 1,
\ee
which yields, together with identity \eqref{two} and assumption \eqref{c1},
$$
0 < \frac{1}{\tau}\left|\frac{P[w]_0-P[w]_{-1}}{w_0-w_{-1}}\right| F(w_0-w_{-1}) \le C
$$
with a constant $C>0$ independent of $x$ and $\tau$. 
For $w_0 = w_{-1}$, we interpret $(P[w]_0-P[w]_{-1})/(w_0-w_{-1})$ as $B_+'(w_{-1})$ or $B_-'(w_{-1})$, according to the notation \eqref{b+}--\eqref{b-}.
Using the fact that $n\tau = T$ we conclude from \eqref{dp3a}--\eqref{dp3c} in both cases that
\be{dp4}
\frac1\tau \sumi (P[w]_{i+1} - 2P[w]_i + P[w]_{i-1})f(w_{i+1} - w_i)
\ge \frac{\beta}{2}\sumi |w_{i+1} - w_i| - C
\ee
with a constant $C>0$ independent of $x$ and $\tau$. 

We now test the difference of \eqref{dp1} taken at discrete times $i+1$ and $i$ by $\vp = f(w_{i+1} - w_i)$ with $f$ from \eqref{dp2}, which yields by virtue of Hypothesis~\ref{hyd}, \eqref{hi}, and \eqref{dp4} that
\begin{align}\nonumber
& \frac{\beta}{2}\sumi\io |w_{i+1} - w_i|\dd x + \sumi\io \nabla(u_{i+1}{-}u_i)\cdot \nabla f(w_{i+1}{-} w_i)\dd x\\ \label{dp0}
&\qquad + \sumi\ipo b(x)(u_{i+1}-u_i)f(w_{i+1} - w_i)\dd s(x) \le C.
\end{align}
From the obvious inequality $(u_{i+1}-u_i)f(w_{i+1} - w_i) \ge 0$, it follows that
\be{dp5}
\frac{\beta}2\sumi \io |w_{i+1} - w_i|\dd x + \sumi\io \nabla(u_{i+1}-u_i)\cdot \nabla f(w_{i+1} - w_i)\dd x \le C
\ee
with a constant $C>0$ independent of $\tau$. We further have
\begin{align*}
&\nabla(u_{i+1}-u_i)\cdot \nabla f(w_{i+1} - w_i) = f'(w_{i+1} - w_i)\nabla(u_{i+1}-u_i)\cdot \nabla (g(u_{i+1}) - g(u_i))\\
&\quad = f'(w_{i+1} - w_i)\Big(g'(u_{i+1})|\nabla(u_{i+1}-u_i)|^2 + (g'(u_{i+1}) - g'(u_i))\nabla(u_{i+1}-u_i)\cdot \nabla u_i \Big).
\end{align*}
By \eqref{hg} and Proposition~\ref{p2} we thus obtain
\begin{align} \nonumber
&\nabla(u_{i+1}-u_i)\cdot \nabla f(w_{i+1} - w_i)\\ \nonumber
&\quad \ge f'(w_{i+1} - w_i)\Big(g_*(\bar U) |\nabla(u_{i+1}-u_i)|^2 - \bar g(\bar U) |u_{i+1} - u_i|\,|\nabla(u_{i+1}-u_i)|\,|\nabla u_i|\Big)\\ \label{dp6}
&\quad \ge \frac{g_*(\bar U)}{2}f'(w_{i+1} - w_i)|\nabla(u_{i+1}-u_i)|^2 - Cf'(w_{i+1} - w_i)|u_{i+1}-u_i|^2|\nabla u_i|^2 
\end{align}
with a constant $C$ independent of $i$ and $\tau$. We have
\be{dp7}
g_*(\bar U)|u_{i+1} - u_i| \le |w_{i+1} - w_i|\le g^*(\bar U)|u_{i+1} - u_i|
\ee
for all $i$, hence
$$
f'(w_{i+1} - w_i) = \frac{\tau}{(\tau + |w_{i+1} - w_i|)^2} \le \frac{\tau}{(\tau + g_*(\bar U) |u_{i+1} - u_i|)^2},
$$
and
$$
\nabla(u_{i+1}-u_i)\cdot \nabla f(w_{i+1} - w_i) \ge \frac{g_*(\bar U)}{2}f'(w_{i+1} - w_i)|\nabla(u_{i+1}-u_i)|^2 - C\tau|\nabla u_i|^2
$$
with a constant $C>0$ independent of $i$ and $\tau$. Combining \eqref{dp5} with \eqref{dp7} and \eqref{es0}, we obtain the assertion.
\epf



The key point in order to perform the actual iteration step is the identification of suitable functions $f_\gamma$ for $\gamma > 1$ in Proposition~\ref{pc}. We first recall the following auxiliary result from the theory of Sobolev spaces, which follows from \eqref{es2b} and the Gagliardo-Nirenberg inequality (see \cite{gagl,nire}).
 
\begin{lemma}\label{dpl1}
Let $\Omega \subset \real^N$ be an open bounded $C^{1,1}$ domain, and put 
\be{ee1}
p_N = 1 + \frac{2}{N}, \quad  \frac1{p_N} + \frac1{p_N'} = 1,\quad \alpha = \frac{\frac12 - \frac 1{2p_N}}{\frac1N} = \frac{N}{N+2} \in (0,1).
\ee
Then there exists $C>0$ such that for all $u\in W^{2,2}(\Omega)$ we have
$$
|\nabla u|_{2p_N} \le C|\nabla u|_2^{1-\alpha} |u|_{2,2,b}^\alpha,
$$
where $|\cdot|_q$ denotes the norm in $L^q(\Omega)$.
\end{lemma}

The next two Lemmas contain the instructions on how to perform the iteration step, with Lemma~\ref{el2} showing how to suitably choose the functions $f_\gamma$ for $\gamma > 1$.

\begin{lemma}\label{el1}
Let the solution to \eqref{dp1} satisfy
$$
\tau^{1-\gamma}\sumi\io |u_{i+1} - u_i|^\gamma\dd x  \le C
$$
for some $\gamma \ge 1$ with a constant $C>0$ independent of $\tau$. Then
$$
\tau^{1-(\gamma/p_N')}\sumi \io |u_{i+1} - u_i|^{\gamma/p_N'}|\nabla u_i|^2(x)\dd x \le C
$$
with $p_N'$ as in \eqref{ee1}, with a constant $C>0$ independent of $\tau$.
\end{lemma}

\bpf{Proof}
By H\"older's inequality we have
$$
\sumi \io |u_{i+1} - u_i|^{\gamma/p_N'}|\nabla u_i|^2\dd x \le \left(\sumi\io |u_{i+1} - u_i|^{\gamma}\dd x \right)^{1/p_N'}\left(\sumi\io |\nabla u_i|^{2p_N}\dd x \right)^{1/p_N}.
$$
By Lemma~\ref{dpl1}, Hypothesis~\ref{hyd} on $u^0 = u_0$, and estimate \eqref{es1} we have
$$
|\nabla u_i|_{2p_N} \le C|\nabla u_i|_2^{1-\alpha} |u_i|_{2,2,b}^\alpha \le C |u_i|_{2,2,b}^\alpha, \ i=0,1,\dots, n-1.
$$
Hence,
$$
\io|\nabla u_i|^{2p_N}\dd x \le C |u_i|_{2,2,b}^{2p_N \alpha} = C |u_i|_{2,2,b}^2,
$$
and by \eqref{es2a} we get
\begin{align*}
&\tau^{1-(\gamma/p_N')}\sumi \io |u_{i+1} - u_i|^{\gamma/p_N'}|\nabla u_i|^2(x)\dd x\\
&\qquad \le C\tau^{1-(\gamma/p_N')}\left(\sumi\io |u_{i+1} - u_i|^{\gamma}\dd x \right)^{1/p_N'}\left(\sumi|u_i|_{2,2,b}^2 \right)^{1/p_N}\\
&\qquad = C\left(\tau^{1-\gamma}\sumi\io |u_{i+1} - u_i|^{\gamma}\dd x \right)^{1/p_N'}\left(\tau\sumi |u_i|_{2,2,b}^2 \right)^{1/p_N} \le C,
\end{align*}
which we wanted to prove.
\epf

\begin{lemma}\label{el2}
Let the solution to \eqref{dp1} satisfy
$$
\tau^{2-\gamma}\sumi \io |u_{i+1} - u_i|^{\gamma-1}|\nabla u_i|^2(x)\dd x \le C
$$
for some $\gamma >1$, with a constant $C>0$ independent of $\tau$. Then
$$
\tau^{1-\gamma}\sumi \io |u_{i+1} - u_i|^{\gamma}\dd x \le C
$$
with a constant $C>0$ independent of $\tau$.
\end{lemma}

\bpf{Proof}
By analogy to \eqref{dp2}, we choose for $\gamma>1$ the function $f$ in Proposition~\ref{pc} as
\be{re7}
f'(w) = f_\gamma'(w) = \tau^{2-\gamma}(\tau+|w|)^{\gamma-3}.
\ee
For $\gamma \ne 2$ we have
$$
\tilde f(w) = wf(w) - F(w) = \int_0^{|w|} v f'(v) \dd v = \tau^{2-\gamma}\int_0^{|w|}(\tau+v)^{\gamma - 2}\dd v - \tau^{3-\gamma}\int_0^{|w|}(\tau+v)^{\gamma - 3}\dd v,     
$$
and by substitution $v = \tau s$ we get
$$
\tilde f(w) =  \tau \int_0^W \left((1+s)^{\gamma-2} - (1+s)^{\gamma-3}\right)\dd s = \tau\int_0^W s\,(1+s)^{\gamma-3} \dd s, \quad W = \frac{|w|}{\tau}.  
$$
For $W \ge 0$ put
$$
\tilde F(W) = W \int_0^W s \,(1+s)^{\gamma-3} \dd s.
$$
We thus have 
$$
\tilde F(W) \ge \frac1{\gamma-1} W^{\gamma} \ \for \ {\gamma} \ge 3.
$$
For $\gamma \in (1,2)\cup (2,3)$ and $W\ge 1$ we estimate
$$
\tilde F(W) = W^{\gamma}\int_0^1 z \,\left(\frac1{W}+z\right)^{\gamma-3} \dd z \ge c_{\gamma} W^{\gamma}, \quad c_{\gamma} = \int_0^1 z \,(1+z)^{\gamma-3} \dd z,
$$
while for $W\in [0,1)$ we have $\tilde F(W) \ge 0 > W^{\gamma} - 1$. We thus have for all $w\in \real$ and $\tau > 0$ that
\be{ckappa}
 |w|(wf(w) - F(w)) \ge \tau^2 C_{\gamma} \big(W^{\gamma} - 1\big)
\ee
with $C_{\gamma} = \min\{1, 1/(\gamma-1), c_{\gamma}\}$. For $\gamma = 2$ we similarly have
$$
wf(w) - F(w) = \int_0^{|w|} v f'(v)\dd v = |w| -
\tau \log\left(1+ \frac{|w|}{\tau}\right) = \tau\big(W - \log(1+W)\big).
$$
Inequality \eqref{ckappa} is satisfied with $C_{\gamma} = C_2 = 1 - \log 2$. Indeed, it is trivial for $W \le 1$, while for $W>1$ it holds with right-hand side $\tau^2 C_2 W^2$. For all $\gamma > 1$ we thus have
\be{tau}
|w|(wf(w) - F(w)) \ge C_{\gamma} \left(\tau^{2-\gamma} |w|^{\gamma} - \tau^2\right).
\ee
We apply again Propositions~\ref{p2} and \ref{pc} with $U = \bar U$ together with \eqref{tau}, and obtain
\begin{align*}
&\frac{1}{\tau}\sumi (P[w]_{i+1} - 2P[w]_i + P[w]_{i-1})f(w_{i+1} - w_i) + \frac{1}{\tau}\frac{P[w]_0 - P[w]_{-1}}{w_0 - w_{-1}}F(w_0 - w_{-1})\\
&\qquad \ge \frac{\beta C_{\gamma}}{2} \left(\Big(\tau^{1-\gamma}\sumi |w_{i+1} - w_i|^{\gamma}\Big) - 1\right).
\end{align*}
Since $f$ is increasing and odd, and $F$ is convex and even, we have for all $w\ne 0$ that $F(w)/|w| = F(|w|)/|w| \le f(|w|)$. Therefore, for $w_0 \neq w_{-1}$ instead of \eqref{dp3b} we now have 
$$
0 < \frac{F(w_0 - w_{-1})}{|w_0 - w_{-1}|} \le \left\{
\begin{array}{ll}
\displaystyle{\frac{1}{2{-}\gamma} \left( 1 - \left(1+ \frac{|w_0 {-} w_{-1}|}{\tau}\right)^{\gamma-2}\right) \le \frac{1}{2-\gamma}} &  \for \gamma<2,\\[2mm]
\displaystyle{\frac{1}{\gamma{-}2} \left(\left(1+ \frac{|w_0 {-} w_{-1}|}{\tau}\right)^{\gamma-2}\! - 1\right) \le \frac{(1+ C g^*(\bar U))^{\gamma-2} - 1}{\gamma-2}} &  \for \gamma>2,\\[2mm]
\displaystyle{\frac{|w_0 - w_{-1}|}{\tau} \le C g^*(\bar U)} &  \for \gamma=2,
\end{array}
\right.
$$
with $C$ from \eqref{umin} and $g^*(\bar U)$ from \eqref{hg}, while for $w_0 = w_{-1}$ we interpret the quotient $(P[w]_0-P[w]_{-1})/(w_0-w_{-1})$ as $B_+'(w_0+)$ or $B_-'(w_0+)$ according to the notation \eqref{b+}--\eqref{b-} and the convention from Lemma~\ref{cl1}. We now proceed as in the proof of Proposition~\ref{pl1}, testing the difference of \eqref{dp1} taken at discrete times $i+1$ and $i$ by $\vp = f(w_{i+1} - w_i)$ with $f$ given by \eqref{re7}. The counterpart of \eqref{dp5} reads
\be{dp5a}
\frac{\beta C_{\gamma}\tau^{{1-\gamma}}}{2}\sumi \io |w_{i+1} - w_i|^{\gamma}\dd x + \sumi\io \nabla(u_{i+1}-u_i)\cdot \nabla f(w_{i+1} - w_i)\dd x \le C
\ee
with a constant $C>0$ independent of $\tau$. We now refer to \eqref{dp6}, where we have by \eqref{dp7} that
$$
f'(w_{i+1} - w_i) = \frac{\tau^{2-\gamma}}{(\tau + |w_{i+1} - w_i|)^{3-\gamma}} \le \frac{\tau^{2-\gamma}}{(\tau + g_*(\bar U) |u_{i+1} - u_i|)^{3-\gamma}},
$$
hence,
\begin{align*}
&\nabla(u_{i+1}-u_i)\cdot \nabla f(w_{i+1} - w_i) \\
&\ge \frac{g_*(\bar U)}{2}f'(w_{i+1} - w_i)|\nabla(u_{i+1}-u_i)|^2 - Cf'(w_{i+1} - w_i)|u_{i+1}-u_i|^2|\nabla u_i|^2
\\
&\ge \frac{g_*(\bar U)}{2}f'(w_{i+1} - w_i)|\nabla(u_{i+1}-u_i)|^2 - C\tau^{2-\gamma} |u_{i+1} - u_i|^{\gamma-1}|\nabla u_i|^2
\end{align*}
with a constant $C>0$ independent of $\tau$ and $i$. Using the above inequality together with \eqref{dp7}, we obtain from \eqref{dp5a} that
\be{dp5b}
\tau^{1-\gamma}\sumi \io |u_{i+1} - u_i|^{\gamma}\dd x \le C \left(1+ \tau^{2-\gamma} \sumi \io |u_{i+1} - u_i|^{\gamma-1}|\nabla u_i|^2 \dd x\right),
\ee
which concludes the proof.
\epf

We are now ready to make the iteration argument explicit. We define a sequence of exponents $\{\gamma_j\}_{j\in\nat\cup\{0\}}$ as
$$
\gamma_0 = 1, \quad \gamma_j := 1 + \frac{\gamma_{j-1}}{p_N'} \ \mbox{ for } j \in \nat,
$$
with $p_N'$ as in \eqref{ee1}. The initial step
$$
\tau^{1-\gamma_0}\sumi\io |u_{i+1} - u_i|^{\gamma_0} \dd x  \le C
$$
corresponds to Proposition~\ref{pl1}. Combining Lemma~\ref{el1} with Lemma~\ref{el2}, for $j \in \nat$ we get the implication
\begin{align*}
\tau^{1-\gamma_j}\sumi\io |u_{i+1} - u_i|^{\gamma_j} \dd x  \le C \ \Longrightarrow \ &\tau^{-\gamma_j/p_N'}\sumi\io |u_{i+1} - u_i|^{1+ (\gamma_j/p_N')}\dd x \\
& = \tau^{1-\gamma_{j+1}}\sumi\io |u_{i+1} - u_i|^{\gamma_{j+1}} \dd x\le C
\end{align*}
with constants $C$ independent of $\tau$, which corresponds to the inductive step. Note that we have
$$
\lim_{j\to \infty} \gamma_j = p_N = 1+ \frac{2}{N}.
$$
For an arbitrary $1\le q < p_N$ we thus obtain, after finitely many steps,
\be{utq}
\tau^{1-q}\sumi\io |u_{i+1} - u_i|^{q}(x)\dd x \le C
\ee
with a constant $C>0$ independent of $\tau$.
We are now ready to prove the existence of a large amplitude solution to \eqref{e0v}.


\section{Proof of Theorem~\ref{t1}}\label{sec:proof}

We construct the sequences $\{\hat u\on\}$ and $\{\bar u\on\}$ of approximations by piecewise linear and piecewise constant interpolations of the solutions to time-discrete problem associated with the division $t_i=i\tau=iT/n$ for $i=0,1,\dots,n$ of the interval $[0,T]$. Namely, for $x \in \Omega$ and $t\in [t_{i-1},t_i)$, $i=1, \dots, n$, we define 
\begin{align}\label{pl}
\hat u\on(x,t) &= u_{i-1}(x) + \frac{t-t_{i-1}}{\tau} (u_i(x) - u_{i-1}(x)),\\ \label{pco}
\bar u\on(x,t) &= u_{i-1}(x),
\end{align}
continuously extended to $t=T$. Similarly, we consider
$$
G\on(x,t) = G[u]_{i-1}(x) + \frac{t-t_{i-1}}{\tau} (G[u]_i(x) - G[u]_{i-1}(x)),
$$
and the data $\bar h\on(x,t) = h(x,t_{i-1}) = h_{i-1}(x)$ for $x \in \Omega$, and ${\bar u}^{*,(n)}(x,t) = u^*(x,t_{i-1}) = u^*_{i-1}(x)$ for $x \in \partial\Omega$, again continuously extended to $t=T$. Then \eqref{de0v} is of the form
\be{de0w}
\io \left(G\on_t\vp + \nabla \bar u\on \cdot\nabla\vp\right)\dd x + \ipo b(x)(\bar u\on - {\bar u}^{*,(n)})\vp \dd s(x) = \io \bar h\on\vp\dd x
\ee
for every test function $\vp \in W^{1,2}(\Omega)$. Our final goal is letting $n \to \infty$ (equivalently, $\tau \to 0$) in \eqref{de0w}, and obtaining in the limit a solution to \eqref{e0v}. This will be achieved by exploiting the a~priori estimates derived in the previous sections.

We start by noting that, by virtue of Lemma~\ref{dpl1} and of \eqref{utq}, the approximate solutions $\hat u\on$ are for each fixed $q\in (1,p_N)$ uniformly bounded in the anisotropic space
\be{aniso}
X_{N,q} := \{u \in L^1(\Omega \times (0,T)): u_t \in L^q(\Omega \times (0,T)), \ \nabla u \in L^{2p_N}(\Omega \times (0,T); \real^N)\}.
\ee

The following compactness result will enable us to pass to the limit as $n\to \infty$.

\begin{proposition}\label{embed}
The sequence $\{\hat u\on; n \in \nat\}$ is compact in the space $L^{r_N}(\Omega; C[0,T])$ for any $r_N \ge 1$ if $N\le 2$, and
$$
1< r_N < r_N^* := \frac{N(N+4)}{N^2 - 4} \ \for N =3, \quad 1< r_N < \hat r_N := \frac{N^2(N+2)}{N^3 - 4(N+2)} \ \for N\ge 4.
$$
\end{proposition}

The proof of Proposition~\ref{embed} is based on the following anisotropic embedding formula as a special case of the general theory in~\cite{bin}.

\begin{proposition}\label{embin}
Let $q_2 \ge \max\{q_0, q_1\}$, and let
\be{emfor}
\left(1 -\frac1{p_0}\right)\left(\frac1N - \frac1{q_1} + \frac1{q_2}\right) > \frac1{p_1}\left(\frac1{q_0} - \frac1{q_2}\right).
\ee
Then the space
$$
Y_N = \left\{u \in L^1(\Omega \times (0,T)): u_t \in L^{q_0}(\Omega; L^{p_0}(0,T)),\, \frac{\partial u}{\partial x_k} \in L^{q_1}(\Omega; L^{p_1}(0,T)) \, \for k=1, \dots, N\right\}
$$
is compactly embedded in $L^{q_2}(\Omega; C[0,T])$.
\end{proposition}

\bpf{Proof of Proposition~\ref{embed}}
Using the estimates in \eqref{aniso}, we apply Proposition~\ref{embin} with
\be{expo1}
q_0 = p_0 = q < p_N,\  p_1 = 2p_N,\ q_1 = \min\{2p_N,r_N\}, \ q_2 = r_N.
\ee
Note that the requirement $q_2 \ge q_1$ is trivially fulfilled. Condition \eqref{emfor} is of the form
\be{expo2}
\left(1 -\frac1{q}\right)\left(\frac1N - \frac{1}{q_1} + \frac{1}{r_N}\right)> \frac1{2p_N}\left(\frac1{q} - \frac1{r_N}\right).
\ee
If $r_N$ and $q_1 = \min\{2p_N,r_N\}$ are such that inequality \eqref{expo2} holds for $q = p_N$, it will certainly hold for $q$ in a sufficiently small left neighborhood of $p_N$. Hence, $r_N$ has to be such that
\be{expo3}
\left(1 -\frac1{p_N}\right)\left(\frac1N - \frac{1}{q_1} + \frac{1}{r_N}\right)> \frac1{2p_N}\left(\frac1{p_N} - \frac1{r_N}\right),
\ee
that is,
\be{expo3a}
\frac{N+4}{r_N} > \frac{N^2}{N+2}- \frac4N + \frac{4}{q_1}.
\ee
Consider first the case $q_1 = 2p_N$. From \eqref{expo3a} it follows that
\be{expo3b}
\frac{N+4}{r_N} > \frac{N^2}{N+2}- \frac4N + \frac{2N}{N+2} = N - \frac4N.
\ee
For $N=1$ or $N=2$, \eqref{expo3b} is satisfied for every $r_N \ge 1$. The condition $q_2 = r_N \ge q_1 = 2p_N$ in combination with \eqref{expo3b} is fulfilled only if
$$
N^3 < 8N + 16,
$$
that is, for $N=3$, and \eqref{expo3b} is equivalent to $r_N < r_N^*$, which is precisely the assertion. For $N\ge 4$ we apply \eqref{expo3a} with $q_1 = r_N$, which is satisfied provided
$$
\frac{1}{r_N} > \frac{N}{N+2} - \frac{4}{N^2},
$$
that is, for $r_N < \hat r_N$. To complete the proof, it suffices to notice that for $N=1,2$ we can always choose $r_N > q$, whereas for every $N \ge 3$ we have $\min\{r_N^*, \hat r_N\} > p_N$. Hence the condition $q_2 = r_N \ge p_N > q_0 = q$ in Proposition~\ref{embin} can always be fulfilled.
\epf

From Proposition~\ref{embed} we conclude that there exists $u \in Y_N$ such that a subsequence of $\{\hat u\on\}$ converges to $u$ weakly in $Y_N$ and strongly in $L^{r_N}(\Omega; C[0,T])$. To prove the convergence of (a subsequence of) $\{\bar u^{(n)}\}$ to the same limit $u$, notice that for $x \in \Omega$ and $t \in [t_{i-1}, t_i)$ we have $|\bar u\on(x,t) - \hat u\on(x,t)| \le |u_i(x) - u_{i-1}(x)|$, hence, 
\be{hatbar}
\int_0^T\!\!\io |\bar u\on - \hat u\on|^q(x,t)\dd x\dd t \le \sumj \io|u_i - u_{i-1}|^q(x) \dd x \le C\tau^{q-1}
\ee
by virtue of \eqref{utq}. We see that $\bar u\on \to u$ strongly in $L^q(\Omega\times (0,T))$ and, by \eqref{es1}, $\nabla \bar u\on \to \nabla u$ weakly-star in $L^\infty(0,T; L^2(\Omega))$. On $\partial\Omega$, we apply the interpolation inequality
$$
\ipo |v|^{q} \dd s(x) \le C \left(\io |v|^{q}\dd x + \io|v|^{q-1}|\nabla v|\dd x\right) \quad \for \ v \in W^{1,2}(\Omega)
$$
to $v = \bar u\on - u$, and conclude that $\bar u\on \to u$ strongly in $L^{q}(\partial\Omega\times (0,T))$ and weakly-star in $L^\infty(0,T; L^2(\partial\Omega))$ by \eqref{es1}. By \eqref{hi0}--\eqref{hi}, the data ${\bar u}^{*,(n)}$, $\bar h\on$ converge strongly to $u^*, h$ in $L^2(\partial\Omega\times (0,T))$ and in $L^2(\Omega\times (0,T))$, respectively. We thus can pass to the limit in the linear terms of \eqref{de0w}, and to check that $u$ is a solution of \eqref{e0v} it remains to prove the convergence $G\on_t \to G[u]_t$.

It follows from \eqref{de5a} and \eqref{es1} that the sequence $\{G\on_t\}$ is bounded in $L^2(\Omega \times (0,T))$. Hence, there exists $G^* \in L^2(\Omega \times (0,T))$ such that a subsequence of $\{G\on_t\}$ converges weakly to $G^*$ as $n \to \infty$. Let us estimate the difference $G\on - G[u]$ in $L^{q}(\Omega; C[0,T])$ for $q \in (1,p_N)$. H\"older's inequality yields
	\begin{equation}\label{convG}
		\begin{aligned}
			&\io \sup_{t\in (0,T)} |G\on(x,t) - G[u](x,t)|^q \dd x \\
			&\le \io \sup_{t\in (0,T)} |G\on(x,t) - G[\hat u\on](x,t)|^q \dd x \\
				&\qquad + |\Omega|^{1-\frac{q}{r_N}} \left(\io \sup_{t\in (0,T)} |G[\hat u\on](x,t) - G[u](x,t)|^{r_N} \dd x\right)^\frac{q}{r_N}.
		\end{aligned}
	\end{equation}
	We estimate the first term at the right-hand side.
Let $v,w:[0,T] \to \real$ be regulated functions, and let the initial memory curve $\lambda$ be the same for $G[v]$ and $G[w]$. Then (see \cite[Proposition~2.3]{jana}) there exists a constant $C_G>0$ such that for all regulated inputs  $v,w : [a,b] \to \real$ we thus have
$$
\sup_{t\in [a,b]}|G[v](t)- G[w](t)| \le C_G \sup_{t\in [a,b]}|v(t)- w(t)|.
$$
Since the initial memory curve $\lambda$ from Definition~\ref{dpr} is the same for $G\on$ and for $G[\hat u]$, we have for all $t \in [0,T]$ and a.\,e. $x \in \Omega$ that
$$
|G\on(x,t) - G[\hat u\on](x,t)| \le C_G \max_{i=1, \dots, n} |u_i(x) - u_{i-1}(x)|.
$$
In particular,
$$
\io \sup_{t\in (0,T)}|G\on(x,t) - G[\hat u\on](x,t)|^{q}\dd x \le C_G\sumj\io|u_i - u_{i-1}|^{q}(x)\dd x \le C\tau^{q-1}
$$
by virtue of \eqref{utq}, with a constant $C>0$ independent of $\tau$. As for the second term at the right-hand side of \eqref{convG}, we employ the Lipschitz continuity of $G$ in $L^{r_N}(\Omega;C[0,T])$ (see \cite[Proposition~II.3.11]{book}), and we conclude that a subsequence of $\{G\on\}$ converges as $n \to \infty$ strongly to $G[u]$ in $L^{q}(\Omega; C[0,T])$. Hence $G^* = G[u]_t$ a.\,e., and the existence proof is complete.

Uniqueness can be proved using the Hilpert inequality already derived in discrete form in \eqref{dh2}, which we now state in the following form.

\begin{proposition}\label{hile}
Let $G$ be a regular Preisach operator with initial memory curve $\lambda$ in the sense of Definition~\ref{dpr}, let $\psi$ be the function defined in \eqref{psi}, and let $H$ be the Heaviside function \eqref{heavi}. Let $u,v \in L^1(\Omega; W^{1,1}(0,T))$ such that $u(x,0) = v(x,0) = \lambda(x,0)$ be arbitrary, and put as in \eqref{e4} $\xi^r = \play_r[\lambda,u]$, $\eta^r = \play_r[\lambda, v]$. Then for almost all $(x,t) \in \Omega\times (0,T)$ we have
$$
(G[u] - G[v])_t H(u-v) \ge \frac{\partial}{\partial t} \int_0^\infty \big(\psi(x,r,\xi^r) - \psi(x,r,\eta^r)\big)^+ \dd r,
$$
where $(\,\cdot\,)^+$ denotes the positive part.
\end{proposition}

\bpf{Proof}
By \eqref{e4a}, the inequality
\be{hi1}
\big(\rho(x,r,\xi^r)\xi^r_t - \rho(x,r,\eta^r)\eta^r_t\big)\,f(u-v) \ge \big(\rho(x,r,\xi^r)\xi^r_t - \rho(x,r,\eta^r)\eta^r_t\big)\,f(\xi^r - \eta^r)
\ee
holds almost everywhere for every nondecreasing function $f:\real \to \real$. In particular, if $f$ is the Heaviside function $H$, then the right-hand side of \eqref{hi1} is of the form
$$
\frac{\partial}{\partial t} \big(\psi(x,r,\xi^r) - \psi(x,r,\eta^r)\big)^+,
$$
and the assertion follows by integrating over $r$ from $0$ to $\infty$.
\epf

We are now ready to finish the proof of Theorem~\ref{t1} by proving that the solution is unique. Indeed, it suffices to test the difference of \eqref{e0v} taken for $u$ and $v$ successively with the regularizations $\vp = H_\ve(u-v)$ and $\vp = H_\ve(v-u)$ of the Heaviside function given by \eqref{hev} for $\ve>0$, and let $\ve$ tend to $0$ to conclude that $u=v$.


\end{document}